\documentclass[amscd,amssymb,verbatim,10pt,letterpaper]{amsart}
\usepackage{latexsym}
\usepackage{euscript}

    \newcommand{\ie}{{\em i.e.}}
    \newcommand{\eg}{{\em e.g.}}

    \newtheorem{thm}{Theorem}[section]
    \newtheorem{prop}[thm]{Proposition}
    \newtheorem{lem}[thm]{Lemma}
    \newtheorem{cor}[thm]{Corollary}

    \theoremstyle{definition}

    \newtheorem{exmp}[thm]{Example}

    \theoremstyle{remark}

    \newtheorem{rem}[thm]{Remark}

    \newcommand{\fn}[3]{#1 \colon #2 \rightarrow #3}

    \newcommand{\fna}[2]{#1 \rightarrow #2}

    \newcommand{\Spec}{\operatorname{Spec}}

     \newcommand{\Hilb}{\operatorname{Hilb}}

    \newcommand{\chark}{\operatorname{char}}

    \newcommand{\script}[1]{\EuScript{#1}}

    \newcommand{\gf}{K}   

\begin{document}

\title[Syzygies of border basis scheme generators]{Some syzygies of the generators of the ideal of a border basis scheme}

\author[M. Huibregtse]{Mark E. Huibregtse}

\address{Department of Mathematics and Computer Science\\
         Skidmore College\\
         Saratoga Springs, New York 12866}

\email{mhuibreg@skidmore.edu}

\date{\today} 


\newcommand{\fmm}{\mathcal{X}}  
\newcommand{\gfmm}{\mathcal{A}} 
\newcommand{\Tr}{\operatorname{Tr}} 
\newcommand{\mdeg}{\operatorname{md}} 
\newcommand{\inc}{\operatorname{inc}} 

\begin{abstract}
     A border basis scheme is an affine scheme that can be viewed as an open subscheme of the Hilbert scheme of $\mu$ points of affine $n$-space.  We study syzygies of the generators of a border basis scheme's defining ideal.  These generators arise as the entries of the commutators of certain matrices (the ``generic multiplication matrices'').  We consider two families of syzygies that are closely connected to these matrices:  The first arises from the Jacobi identity, and the second from the fact that the trace of a commutator is $0$.  Several examples of both types of syzygy are presented, including a proof that the border basis schemes in case $n$ $=$ $2$ are complete intersections.  
\end{abstract}

\maketitle

\section{Introduction} \label{sec:Intro}

	Let $\gf$ be a field, and $\script{O}$ an \textbf{order ideal};
 that is, $\script{O}$ is a finite set
\[
    \script{O} = \{t_1, t_2, \dots, t_{\mu} \}\ \subseteq\ \mathbb{T}^n = \{x_1^{\alpha_1}x_2^{\alpha_2}\dots x_n^{\alpha_n} \mid \alpha_i \geq 0  \}
\]
such that if $m$ $\in$ $\mathbb{T}^n$ and $m \mid t_i$ for some $1$ $\leq$ $i$ $\leq$ $\mu$, then $m$ $\in$ $\script{O}$; we view $\script{O}$ as a subset of the polynomial ring $\gf[x_1, x_2, \dots, x_n]$ $=$ $\gf[\mathbf{x}]$ in the obvious way.  The $\script{O}$\textbf{-border basis scheme} $\mathbb{B}_{\script{O}}$ is an affine scheme that parameterizes the ideals $I$ $\subseteq$ $\gf[\mathbf{x}]$ such that the quotient $\gf[\mathbf{x}]/I$ is $\gf$-free with basis $\script{O}$; as such, it is an open subscheme of the Hilbert scheme of $\mu$ points of $\mathbb{A}^n$.  It is defined as follows: One first computes the \textbf{border} 
\[
    \partial \script{O}\ =\ (x_1\script{O}\cup \dots \cup x_n \script{O})\setminus \script{O}\ =\ \{b_1, \dots, b_{\nu}\}
\]
of the order ideal, and then observes that if the quotient $\gf[\mathbf{x}]/I$ is $\gf$-free with basis $\script{O}$, then every monomial $b_j$ $\in$ $\partial \script{O}$ is congruent to a unique $\gf$-linear combination of the monomials $t_i$ $\in$ $\script{O}$; that is, one has unique polynomials $f_j$ $\in$ $I$ of the form
\[
    f_j = b_j - \sum_{i=1}^{\mu} a_{ij}\,t_i \in I, \ \ a_{ij} \in \gf,
\]
that in fact form a basis of $I$ called the $\script{O}$\textbf{-border basis}.  One now replaces the scalars $a_{ij}$ by indeterminates $c_{ij}$, and defines polynomials
\[
    g_j = b_j - \sum_{i=1}^{\mu} c_{ij}\,t_i \in \gf[c_{11}, \dots, c_{\nu \mu}][\mathbf{x}] = \gf[\mathbf{c}][\mathbf{x}];
\]
the ``generic quotient with basis $\script{O}$'' is then
\begin{equation} \label{eqn:genquot}
    \gf[\mathbf{c}][\mathbf{x}]/(g_{ij}).
\end{equation}
One then computes the $(\mu \times \mu)$-matrices $\gfmm_k$ of the linear maps induced by multiplication by $x_k$ on the quotient (\ref{eqn:genquot}) with respect to the ``basis'' $\script{O}$; that is, one computes formally, assuming that $\script{O}$ is a $\gf[\mathbf{c}]$-basis.  If one now replaces the indeterminates $c_{ij}$ with scalars $a'_{ij}$, one obtains an ideal 
\[
    I' = (g_j(a'_{ij})) \subseteq \gf[\mathbf{x}],
\]
and it is known (see, \eg, \cite[Th.\ 6.4.30, p.\ 434]{KreutzerAndRobbianoVolTwo}) that the quotient $\gf[\mathbf{x}]/I'$ is $\gf$-free with basis $\script{O}$ if and only if the various multiplication maps $\gfmm_{k}(a'_{ij})$ on the quotient commute with one another.  Therefore, the entries of the commutators 
\begin{equation} \label{eqn:rhoDefn}
    [\gfmm_k, \gfmm_l]\ =\ (\rho^{kl}_{pq}),\ \  1 \leq k < l \leq n,\ 1 \leq p, q \leq \mu,
\end{equation}
which are polynomials in the indeterminates $c_{ij}$, give the polynomial relations that must be satisfied by the tuple of scalars $(a'_{ij})$ in order that the ideal $I'$ yield a $\gf$-free quotient with basis $\script{O}$.  Following \cite[Sec.\ 3]{KreutzerAndRobbiano1:DefsOfBorderBases}, one defines the $\script{O}$-border basis scheme to be
\[
    \mathbb{B}_{\script{O}}\ =\ \Spec(\gf[\mathbf{c}]/I(\mathbb{B}_{\script{O}}),\ \ I(\mathbb{B}_{\script{O}}) = ( \rho^{kl}_{pq} ).
\]

This paper presents two families of syzygies among the generators $\rho^{kl}_{pq}$ of the ideal $I(\mathbb{B}_{\script{O}})$.  These arise naturally from the definition of the generators as entries of the commutators $[\gfmm_k, \gfmm_l]$.  The first family, defined and studied in Section \ref{sec:JISyz}, comes from the well-known Jacobi identity, which implies that (for $k$ $<$ $l$ $<$ $m$)
\[
  \begin{array}{c}
    [\gfmm_k, [\gfmm_l, \gfmm_m]] + [\gfmm_l, [\gfmm_m, \gfmm_k]] + [\gfmm_m, [\gfmm_k, \gfmm_l]]\ =\ 0,\ \text{or} \vspace{.05in}\\
 {[}\gfmm_k, (\rho^{lm}_{pq})]] - [\gfmm_l, (\rho^{km}_{pq})] + [\gfmm_m, (\rho^{kl}_{pq})]\ =\ 0;
  \end{array}
\]
whence, every entry on the LHS is a (possibly trivial) syzygy of the $\rho$'s.  Several examples of these \textbf{Jacobi identity syzygies} are presented in Section \ref{sec:JISyzExs}.  (Note that there are no nontrivial Jacobi identity syzygies in two dimensions, that is, when one is working with order ideals $\script{O}$ $\subseteq$ $\mathbb{T}^2$.)

The second family of syzygies, defined and studied in Section \ref{sec:Tsyz}, following a brief technical preparation in Section \ref{sec:sumOfCom}, is based on the fact that the trace of a commutator is $0$.  In particular, one has that (for $k$ $<$ $l$)
\[
    \Tr([\gfmm_k, \gfmm_l])\ =\ \Tr(\rho^{kl}_{pq})\ =\ \rho^{kl}_{1,1} + \rho^{kl}_{2,2} + \dots + \rho^{kl}_{\mu, \mu}\ =\ 0;
\]
these are the simplest of the \textbf{trace syzygies}.  Others are obtained by forming expressions that are sums of products of the $\gfmm_k$ and the $[\gfmm_k, \gfmm_l]$ that reduce to commutators, and therefore have trace $0$; these expressions are presented in Section \ref{sec:sumOfCom}.  As a simple example, one has that (for $k$ $<$ $l$)
\[
   \Tr\left([\gfmm_k, \gfmm_l]\, \gfmm_l + \gfmm_l\,[\gfmm_k, \gfmm_l]\right)\ =\ 0,
\]
since the matrix expression reduces to the commutator $[\gfmm_k, \gfmm_l\, \gfmm_l]$; moreover, the entries in the matrix are again clearly $\gf[\mathbf{c}]$-linear combinations of the $\rho^{kl}_{pq}$ that equal $0$, that is, syzygies of the $\rho$'s.

The final three sections of the paper present extended examples of the trace syzygies.  The first of these, in Section \ref{sec:TSyzEx1}, studies the trace syzygies when $\script{O}$ $=$ $\{1, x_1, x_2 \}$ $\subseteq$ $\mathbb{T}^2$; certain syzygies first presented in \cite[Sec.\ 5.2]{Huib:HaimanUmu}, and sufficient to demonstrate that $\mathbb{B}_{\script{O}}$ is isomorphic to six-dimensional affine space, are re-derived using trace syzygies. (The derivation of the syzygies in the previous paper is far less general than that given here.)  

Section \ref{sec:TSyzEx2} computes various trace syzygies in case $\script{O}$ $=$ $\{ 1, x_1 \}$ $\subseteq$ $\mathbb{T}^3$; here again, one has that $\mathbb{B}_{\script{O}}$ is isomorphic to a six-dimensional affine space, but the demonstration requires both trace syzygies and Jacobi identity syzygies (the latter are worked out for this case in Example \ref{exmp:JacId3}).  It should be noted that the result of Section  \ref{sec:TSyzEx2} is (after renumbering the variables) a special case of \cite[Cor.\ 3.13]{Robbiano:BorderAndGrobner}, which states that whenever 
$\script{O}$ $=$ $\{1,\, x_n,\, x_n^2,\, \dots,\, x_n^{\mu-1} \}$,
then $\mathbb{B}_{\script{O}}$ is isomorphic to the affine space $\mathbb{A}^{\mu n}$.

The final Section \ref{sec:2varGenEx} considers the trace syzygies for an arbitrary order ideal $\script{O}$ $\subseteq$ $\mathbb{T}^2$; in case $\chark(K)$ $=$ $0$, it is shown  that $\mathbb{B}_{\script{O}}$ is a complete intersection in $\Spec(\gf[\mathbf{c}])$ $=$ $\mathbb{A}^{\mu \nu}$, and a minimal set of generators of the ideal $I(\mathbb{B}_{\script{O}})$ is specified.

The sections not yet mentioned concern various preliminaries: Section \ref{sec:BorderBasisSchemes} recalls the definition of border basis schemes in somewhat more detail than above; Section \ref{sec:BBSidGens} takes a closer look at the polynomials $\rho^{kl}_{pq}$ that generate the ideal $I(\mathbb{B}_{\script{O}})$; and Section \ref{sec:cDegs} discusses an assignment of multi-degrees to the indeterminates $c_{ij}$ under which our syzygies are homogeneous. 

Before commencing, I am pleased to acknowledge my debt to the beautiful articles \cite{Haiman:CN-HS}, \cite{KreutzerAndRobbiano1:DefsOfBorderBases} and \cite{Robbiano:BorderAndGrobner}.  I also thank the referee for several helpful comments that greatly improved the exposition.

\section{Border basis schemes} \label{sec:BorderBasisSchemes}

	Border basis schemes are defined in \cite[Sec.\ 3]{KreutzerAndRobbiano1:DefsOfBorderBases}; we repeat their definition here,  using the same notation.  Recall that an \textbf{order ideal} is a finite set of monomials 
\[
  \script{O} = \{t_1, t_2, \dots, t_{\mu} \}\ \subseteq\ \mathbb{T}^n = \{x_1^{\alpha_1}x_2^{\alpha_2}\dots x_n^{\alpha_n} \mid \alpha_i \geq 0  \}]
\] such that $m$ $\in$ $\mathbb{T}^n$ and $m \mid t_i$ for some $1$ $\leq$ $i$ $\leq$ $\mu$  implies $m$ $\in$ $\script{O}$, and the \textbf{border} of $\script{O}$ is 
\[
   \partial \script{O}\ =\ (x_1\script{O}\cup \dots \cup x_n \script{O})\setminus \script{O}\ =\ \{b_1, \dots, b_{\nu}\}.
\]
 The $\script{O}$-border basis scheme parameterizes the ideals $I$ $\subseteq$ $\gf[\mathbf{x}]$ such that the quotient $\gf[\mathbf{x}]/I$ has $\script{O}$ as $\gf$-basis.  Such an $I$ must have a (unique) basis of the form $\{f_1,\dots,f_{\nu}\}$, where 
\[
  f_j = b_j - \sum_{i=1}^{\mu}a_{ij}t_i,\ \text{with } a_{ij} \in \gf.
\]
Of course, the coefficients $a_{ij}$ cannot be chosen arbitrarily; in order for the quotient to be $\gf$-free with basis $\script{O}$, it is necessary and sufficient that the the matrices $\fmm_k$ representing the $\gf$-linear maps defined by multiplication by the variables $x_k$ on $\gf[\mathbf{x}]/I$ with respect to the basis $\script{O}$ are pairwise commutative (see, \eg,  \cite[Th.\ 6.4.30, p.\ 434]{KreutzerAndRobbianoVolTwo}).  It follows that the tuples $(a_{ij})$ such that the corresponding quotient $\gf[\mathbf{x}]/I$ is $\gf$-free with basis $\script{O}$ are the points of an affine algebraic set.

	To make this explicit, we introduce a set of indeterminates 
\[
    \{c_{ij} \mid\ 1 \leq i \leq \mu,\ 1 \leq j \leq \nu \},
\]
and define polynomials
\begin{equation} \label{eqn:gpolys}
    g_j = ( b_j - \sum_{i=1}^{\mu} c_{ij}t_i) \in \gf[c_{11},\dots,c_{\mu \nu}, \mathbf{x}] = \gf[\mathbf{c}, \mathbf{x}], \ 1 \leq j \leq \nu;
\end{equation}
the set $\{ g_1, \dots, g_{\nu} \}$ is called the \textbf{generic $\script{O}$-border prebasis} in \cite[Def.\ 3.1]{KreutzerAndRobbiano1:DefsOfBorderBases}.  We then form the generic multiplication matrices $\gfmm_k$ that formally correspond to multiplication by the variables $x_k$ on the quotient $\gf[\mathbf{c}, \mathbf{x}]/(g_j)$ with respect to the ``basis'' $\script{O}$ (that is, one computes as if $\script{O}$ really were a basis of the quotient). We have that $\gfmm_k$ $=$ $(\xi^{k}_{rs})$ is an $n\times n$-matrix with entries $\xi^{k}_{rs}$ $\in$ $\gf[\mathbf{c}]$;
one sees easily that 
\begin{equation} \label{eqn:gfmmDescrip} 
    \xi^{k}_{rs} = \left\{
     \begin{array}{l}
	0, \textrm{ if } x_k\cdot t_s = t_{s'} \in \script{O}, \textrm{ and } s' \neq r;\vspace{.05in} \\
        1, \textrm{ if } x_k\cdot t_s = t_{s'} \in \script{O}, \textrm{ and } s' = r;\vspace{.05in}\\
        c_{rj}, \textrm{ if } x_k\cdot t_s = b_j \in \partial \script{O}.
     \end{array} \right.
\end{equation}
The \textbf{$\script{O}$-border basis scheme} is then ``[t]he affine scheme $\mathbb{B}_{\script{O}}$ $\subseteq$ $\mathbb{A}^{\mu \nu}$ defined by the ideal $I(\mathbb{B}_{\script{O}})$ generated by the entries of the matrices $\gfmm_k \gfmm_l - \gfmm_l \gfmm_k$ with $1$ $\leq$ $k$ $<$ $l$ $\leq$ $n$\dots'' \cite[Def.\ 3.1]{KreutzerAndRobbiano1:DefsOfBorderBases}.  Alternate constructions of the border basis schemes are given in \cite{Huib:UConstr} and \cite{LaksovAndSkjelnes:HilbSchConstr}.

\section{Ideal generators for the $\script{O}$-border basis scheme} \label{sec:BBSidGens}

We have just seen that the ideal of the $\script{O}$-border basis scheme is generated by the entries of the commutators 
\[
    [\gfmm_k, \gfmm_l]\ =\ \gfmm_k \gfmm_l - \gfmm_l \gfmm_k,\ 1 \leq k < l \leq n;
\]
we denote the $(p,q)$-entry of this matrix by 
\[
    \rho^{kl}_{pq} \in \gf[\mathbf{c}], 1 \leq p, q \leq \mu.
\]
These entries are clearly polynomials of total degree $\leq 2$ in the indeterminates $c_{ij}$.  


  We proceed to write down explicit expressions for the polynomials $
\rho^{kl}_{pq}$.  To begin, one has that
\begin{equation} \label{eqn:rhoForm1}
  \rho^{kl}_{pq} = (\gfmm_k \gfmm_l - \gfmm_l \gfmm_k)_{(p,q)} = \sum_{i=1}^{\mu}\left( \xi_{pi}^{k}\,\xi_{iq}^l - \xi_{pi}^l\,\xi_{iq}^k \right). 
\end{equation}
We want to rewrite this expression in terms of the indeterminates $c_{ij}$; to this end, it is convenient to define the following functions $\sigma_k$, $\tau_k$, $1$ $\leq$ $k$ $\leq$ $n$, and their ``inverses'' $\sigma_k^{\prime}$, $\tau_k^{\prime}$, that encode certain aspects of multiplication by the variables $x_k$ on $\gf[\mathbf{x}]$:
\begin{equation} \label{eqn:auxfuncs}
  \begin{array}{cc}
	\fn{\sigma_k}{\{1,2\dots,\mu \}}{\{0,1,2,\dots,\nu \}} & 
        \fn{\sigma_k^{\prime}}{\{1, 2, \dots, \nu  \}}{\{0,1,\dots, \mu \}}\vspace{.05in}\\
	\sigma_k(i) = \left\{ \begin{array}{l}
                                 j,\ \text{if } x_k \cdot t_i = b_j \in \partial \script{O},\\
0,\ \text{otherwise},
                              \end{array}\right. & \sigma_k^{\prime}(j) = 
     \left\{ \begin{array}{l}
          i, \text{ if } b_j/x_k = t_i \in \script{O},\\
          0, \text{ otherwise},
     \end{array} \right. ,
\vspace{.10in}\\
	\fn{\tau_k}{\{1,2\dots,\mu \}}{\{0, 1,2,\dots,\mu \}} & \fn{\tau_k^{\prime}}{\{ 1,2,\dots,\mu \}}{\{0,1,\dots,\mu\}}\vspace{.05in}\\
\tau_k(i) = \left\{ \begin{array}{l}  
                       i_1,\ \text{if } x_k \cdot t_i = t_{i_1} \in \script{O},\\
                       0,\ \text{otherwise},
                    \end{array}
            \right.	           
                      & \tau_k^{\prime}(i_1) = \left\{ \begin{array}{l}
                                  i,\ \text{if } t_{i_1}/x_k = t_{i}\ \in                                       \script{O},\\
                                  0,\ \text{otherwise}.
                                                     \end{array}\right.
  \end{array}
\end{equation}


It is also convenient to adopt the following conventions, \emph{which we shall do henceforth}:
\begin{equation}\label{eqn:conv}
  \begin{array}{c}
   a_{i0} = a_{0j}= c_{i0} = c_{0j} = 0,\ 1 \leq i \leq \mu,\ 1 \leq j \leq \nu,\vspace{.05in}\\
   b_0 = 0,\ f_0 = (b_0 - \sum_{i=1}^{\mu}a_{i0}t_i) = g_0 = (b_0 - \sum_{i=1}^{\mu}c_{i0}t_i) = 0,\ \text{and}\vspace{.05in}\\ 
  \rho^{kl}_{0q} =  \rho^{kl}_{p0} = 0,\ 0 \leq k < l \leq n,\ 0 \leq p, q \leq \mu.
  \end{array}
\end{equation}
We may now proceed to rewrite the ideal generators $\rho^{kl}_{pq}$ (\ref{eqn:rhoForm1}) using (\ref{eqn:gfmmDescrip}), (\ref{eqn:auxfuncs}), and (\ref{eqn:conv}); 
as explained in the proof of \cite[Prop.\ 6.4.32, p.\ 436]{KreutzerAndRobbianoVolTwo}, there are four cases to consider:
\begin{description}
  \item[Case 1] $x_k\ \cdot t_q$ $=$ $t_{r}$, $x_l \cdot t_q$ $=$ $t_{s}$, $x_k \cdot t_{s}$ $=$ $x_l \cdot t_{r}$ $=$ $t_p \in \script{O}$. One has   
    \[
      \rho^{kl}_{pq} = \xi^k_{p, s}\, \xi^l_{s,q} - \xi^l_{p,r} \xi^k_{r, q} = 1 \cdot 1 - 1 \cdot 1= 0.
    \]   
  \item[Case 2] $x_k\ \cdot t_q$ $=$ $t_{r}$, $x_l \cdot t_q$ $=$ $t_{s}$, $x_k \cdot t_{s}$ $=$ $x_l \cdot t_{r}$ $=$ $b_j$ $\in \partial \script{O}$. One has   
    \[ \rho^{kl}_{pq} = \xi^k_{p,s} \xi^l_{s,q} - \xi^l_{p,r} \xi^k_{r,q} = c_{pj} \cdot 1 - c_{pj} \cdot 1 = 0.
    \]
  \item[Case 3] $x_k \cdot t_q = t_r \in \script{O}$, $x_l \cdot t_q = b_{j_1} \in \partial \script{O}$. In this case, note that 
\[
   x_l\,t_r = x_l\,x_k\,t_q = x_k\,x_l\,t_q = x_k\,b_{j_1} = b_{j_2} \in \partial \script{O}.  
\]
Therefore (keeping in mind our conventions (\ref{eqn:conv})), 
\begin{equation}\label{eqn:Case3}
  \begin{array}{c}
 \rho^{kl}_{pq} = c_{\tau_k^{\prime}(p),j_1} + \sum_{i=1}^{\mu} c_{p,\sigma_k(i)}\cdot c_{i,j_1} - c_{p,j_2},\vspace{.05in}\\
 \text{where } j_1 = \sigma_l(q),\ j_2 = \sigma_l(\tau_k(q)).
  \end{array}.
\end{equation}
\item[Case 4] $x_l \cdot t_q = b_{j_1} \in \partial \script{O}$, $x_k \cdot t_q = b_{j_2} \in \partial \script{O}$.  In this case both the positive and negative terms are computed in the same way as are the positive terms in Case 3; we obtain
\[ 
  \begin{array}{c}
    \rho^{kl}_{pq} = \left( c_{\tau_k^{\prime}(p),j_1} + \sum_{i=1}^{\mu} c_{p,\sigma_k(i)}\cdot c_{i,j_1}  \right) - \left( c_{\tau_l^{\prime}(p),j_2} + \sum_{i=1}^{\mu} c_{p,\sigma_l(i)}\cdot c_{i,j_2} \right),\vspace{.05in}\\
    \text{where } j_1 = \sigma_l(q),\ j_2 = \sigma_k(q).
  \end{array}
\]
\end{description}

In particular, the polynomial $\rho^{kl}_{pq}$ $=$ $0$ if it falls in Case 1 or Case 2 above; in these cases, we say that $\rho^{kl}_{pq}$ is \textbf{trivially} $\mathbf{0}$.  The following example shows that it is possible for $\rho^{kl}_{pq}$ to reduce to $0$ even though it is not, as defined,  trivially $0$:

\begin{exmp} \label{exmp:zeroRhos}
    Let $\script{O}$ $=$ $\{t_1$ $=$ $1\}$ $\subseteq$ $\mathbb{T}^n$.  Then $\partial \script{O}$ $=$ $\{x_1, \dots, x_n\}$, and $\gfmm_k$ $=$ $(c_{1,k})$ for $1$ $\leq$ $k$ $\leq$ $n$.  For $k$ $<$ $l$, we have that $[\gfmm_k,\gfmm_l]$ $=$ $(\rho^{kl}_{1,1})$, and $\rho^{kl}_{1,1}$ is not trivially $0$ (as we are in Case 4); however, we have that 
\[
    (\rho^{kl}_{1,1})\ =\ [\gfmm_k, \gfmm_l]\ =\ (c_{1,k})(c_{1,l}) - (c_{1,l})(c_{1,k})\ = (0).
\]
Thus, we see that $\mathbf{B}_{\script{O}}$ $=$ $\Spec(\gf[c_{1,1}, c_{1,2}, \dots, c_{1,n}]/(0))$ $=$ $\mathbb{A}^n$ in this (simplest) case.
\end{exmp}

\begin{cor} \label{cor:cDegs}
\noindent
\begin{enumerate}
  \item The non-zero terms (if any) in the polynomial $\rho^{kl}_{pq}$ have either total degree $1$ or total degree $2$ in the variables $c_{ij}$; moreover, there are at most two degree-$1$ terms. 
  \item We have that $\rho^{kl}_{pq}$ is trivially $0$ if and only if $x_k t_q$ $\in$ $\script{O}$ and $x_l t_q$ $\in$ $\script{O}$.  \qedsymbol 
\end{enumerate}
\end{cor} 

\newcommand{\Zc}{\mathbb{Z}[\mathbf{c}]}
\begin{rem} \label{rem:rhosInZc}
    The ideal generators $\rho^{kl}_{pq}$ can be viewed as elements of $\Zc$, and most of the work in this paper will be done over $\mathbb{Z}$.  In particular, the Jacobi identity and Trace syzygies will be defined as $\Zc$-syzygies of the $\rho$'s.
\end{rem}

\section{A multi-graded ring structure for $\mathbb{Z}[\mathbf{c}, \mathbf{x}]$}\label{sec:cDegs}
 
We begin by making $\mathbb{Z}[\mathbf{x}]$ into a multi-graded ring in the usual way, by assigning the variable $x_k$ the multi-degree 
\[
    \mdeg(x_k) = e_k^{tr} = (0,\dots, 0, 1, 0,\dots, 0)^{tr},
\]
the $k^{\text{th}}$ column of the $n \times n$ identity matrix; 
the monomial 
$m$ $=$ $x_1^{r_1}\,x_2^{r_2}\dots x_n^{r_n}$ $\in$ $\mathbb{Z}[\mathbf{x}]$ then has multi-degree 
\[
    \mdeg(m)  = (r_1, r_2, \dots, r_n)^{tr}.
\]
We now assign multi-degrees to the indeterminates $c_{ij}$ so that each  polynomial $g_j$ $=$ $b_j - \sum_{i=1}^{\mu}c_{ij}t_i$ in the generic $\script{O}$-border prebasis (\ref{eqn:gpolys}) becomes homogeneous of multi-degree $\mdeg(b_j)$; that is, we define
\[
   \mdeg(c_{ij})\ = \ \mdeg(b_j)  - \mdeg(t_i).
\]
Note that if $\gamma$ $\in$ $\mathbb{Z}[\mathbf{c}, \mathbf{x}]$ is homogeneous of multi-degree $\mathbf{r}$ $=$ $(r_1, \dots, r_n)^{tr}$, then $x_k \cdot \gamma$ is homogeneous of multi-degree 
\[
   \inc_k(\mathbf{r}) := \mathbf{r} + e_k^{tr} = (r_1, \dots, r_{k-1}, r_k + 1, r_{k+1}, \dots, r_n)^{tr}.
\]

Let $F$ denote the free graded $\mathbb{Z}[\mathbf{c}, \mathbf{x}]$-module with basis $\script{O}$ $=$ $\{t_1, t_2, \dots, t_{\mu}\}$,
and endow each basis element with the multi-degree $\mdeg(t_i)$ assigned above.  We interpret the generic multiplication matrix $\gfmm_k$ $=$ $(\xi^k_{pq})$ as the matrix of a $\mathbb{Z}[\mathbf{c}, \mathbf{x}]$-linear map $\fna{F}{F}$ with respect to the basis $\script{O}$.  It is straightforward to check that $\xi^k_{pq}$ (recall (\ref{eqn:gfmmDescrip})) is homogeneous of multi-degree
\[
    \mdeg(\xi^k_{pq}) = \inc_k(\mdeg(t_q)) - \mdeg(t_p);
\]
in particular, $\gfmm_k$ is homogeneous (see \cite[4.7.1]{KreutzerAndRobbianoVolTwo}); consequently, by \cite[4.7.4]{KreutzerAndRobbianoVolTwo}, the linear map defined by $\gfmm_k$ is homogeneous (that is, it maps homogeneous elements to homogeneous elements).  In this case, one sees that if $v$ $\in$ $F$ is homogeneous of multi-degree $\mdeg(v)$ $=$ $\mathbf{r}$, then $\gfmm_k \cdot v$ is homogeneous of multi-degree $\inc_k(\mathbf{r})$.  As an immediate consequence, we obtain 

\begin{lem} \label{lem:rhoMultDeg}
Let $\mathbf{\alpha}$ $=$ $(\alpha_1, \dots, \alpha_{\mu})^{tr}$ be a column vector of elements $\alpha_i$ $\in$ $\mathbb{Z}[\mathbf{c}]$ such that the associated polynomial $\sum_{i=1}^{\mu}\alpha_i t_i$ is homogeneous of multi-degree 
 $\mathbf{r}$ $=$ $(r_1, \dots, r_n)^{tr}$.  Let $\beta$ $=$ $(\beta_1, \dots, \beta_{\mu})^{tr}$ denote the column vector $\gfmm_k \cdot \alpha$.  Then
\begin{enumerate}
  \item The polynomial $\sum_{i=1}^{\mu}\beta_i t_i$ is homogeneous of multi-degree $\inc_k(\mathbf{r})$.
  \item The generators $\rho^{kl}_{pq}$ of the ideal of the $\script{O}$-border basis scheme are homogeneous of multi-degree
\[
    \mdeg(\rho^{kl}_{pq})\ = \ \inc_k(\inc_l(\mdeg(t_q))) - \mdeg(t_p).
\]
\end{enumerate}
\qedsymbol
\end{lem}

\section{The Jacobi identity syzygies} \label{sec:JISyz}

With this section, we reach our main subject, which is the study of certain syzygies among the polynomials $\rho^{kl}_{pq}$.  Here we consider our first source of syzygies, the Jacobi identity.  Later sections will take up our second source of syzygies, which, briefly stated, is the fact that the trace of a commutator is 0.

Let $M_1$, $M_2$, and $M_3$ be square matrices over a commutative ring.  The \textbf{Jacobi identity} is the well-known equality
\[
  [M_1, [M_2, M_3]] + [M_2, [M_3, M_1]] + [M_3, [M_1, M_2]]\ =\ 0,
\]
where, as usual, $[M,N]$ denotes the commutator $M \cdot N - N \cdot M$.

In particular, if we select three variables $x_k$, $x_l$, $x_m$ and apply the Jacobi identity to the associated multiplication matrices $\gfmm_k$, $\gfmm_l$, $\gfmm_m$, we obtain 
\[
        [\gfmm_k, [\gfmm_l, \gfmm_m]] + [\gfmm_l, [\gfmm_m, \gfmm_k]] + [\gfmm_m, [\gfmm_k, \gfmm_l]]   = 0. 
\]
Assuming that the variables are listed in increasing order of index $k$ $<$ $l$ $<$ $m$, and recalling (\ref{eqn:rhoDefn}), we can rewrite this equality as
\[
       [\gfmm_k, (\rho^{lm}_{pq})] - [\gfmm_l, (\rho^{km}_{pq})] + [\gfmm_m, (\rho^{kl}_{pq})]  =  0;
\]
recalling  (\ref{eqn:gfmmDescrip}) and Remark \ref{rem:rhosInZc}, we note that the last equation holds in the ring of $(\mu \times \mu)$-matrices over $\Zc$.
 Therefore, for $1$ $\leq$ $p, q$ $\leq$ $\mu$, the $(p,q)$-entry on the left of the last equality is a $\Zc$-linear combination of the $\rho$'s that is equal to $0$ $\in$ $\Zc$; that is, a $\Zc$-syzygy of the $\rho$'s.  Writing this entry out explicitly, we have that

\[
    \sum_{i=1}^{\mu}\left( (\xi^k_{pi}\rho^{lm}_{iq} - \rho^{lm}_{pi}\xi^k_{iq}) - 
    (\xi^l_{pi}\rho^{km}_{iq} - \rho^{km}_{pi}\xi^l_{iq}) +
    (\xi^m_{pi}\rho^{kl}_{iq} - \rho^{kl}_{pi}\xi^m_{iq}) \right)\ =\ 0,
\]
which can be rewritten, using (\ref{eqn:gfmmDescrip}), (\ref{eqn:auxfuncs}), and the conventions (\ref{eqn:conv}), as 
\begin{equation} \label{eqn:Jsyzexp2}
    \begin{array}{lcl}
       (\rho^{lm}_{\tau_k^{\prime}(p),q} + \sum_{i=1}^{\mu}c_{p,\sigma_k(i)} \rho^{lm}_{iq}) -(\rho^{lm}_{p,\tau_k(q)} + \sum_{i=1}^{\mu}c_{i,\sigma_k(q)}\rho^{lm}_{pi})\ - & {} & {}\vspace{.05in}\\
       (\rho^{km}_{\tau_l^{\prime}(p),q} + \sum_{i=1}^{\mu}c_{p,\sigma_l(i)} \rho^{km}_{iq}) +(\rho^{km}_{p,\tau_l(q)} + \sum_{i=1}^{\mu}c_{i,\sigma_l(q)}\rho^{km}_{pi})\ + & {} & {}\vspace{.05in}\\
       (\rho^{kl}_{\tau_m^{\prime}(p),q} + \sum_{i=1}^{\mu}c_{p,\sigma_m(i)} \rho^{kl}_{iq}) -(\rho^{kl}_{p,\tau_m(q)} + \sum_{i=1}^{\mu}c_{i,\sigma_m(q)}\rho^{kl}_{pi}) & {=} & {0}.
    \end{array}
\end{equation}

Of course, some of the $\rho$'s in the preceding equation might be trivially $0$, and so ignorable (recall from Corollary \ref{cor:cDegs} that $\rho^{kl}_{pq}$ is trivially $0$ if $x_k t_q \in \script{O}$ and $x_l t_q \in \script{O}$).  Hence, once and for all, we make a complete list of the $\rho$'s that are \emph{not} trivially $0$ in some order 
$\rho_1$, $\rho_2$, \dots, $\rho_{\omega}$,
and define $J^{klm}_{pq}$ to be the corresponding tuple 
\begin{equation}\label{eqn:JsyzDef}
  J^{klm}_{pq}\ =\   (\kappa_1, \kappa_2, \dots, \kappa_{\omega}),\ 1 \leq k < l < m \leq n,\ 1 \leq p, q \leq \mu,
\end{equation}
where $\kappa_{\iota}$ $\in$ $\Zc$ is the coefficient of $\rho_{\iota}$ on the LHS of (\ref{eqn:Jsyzexp2}).  We call these syzygies $J^{klm}_{pq}$ the \textbf{Jacobi identity syzygies}.  

\begin{rem} \label{rem:noJacIn2Vars} There are no Jacobi identity syzygies in case $n$ $\leq$ $2$, since their definition requires three distinct variables $x_k$, $x_l$, $x_m$.
\end{rem}

More generally, the syzygies of the $\rho$'s considered in this paper are  defined by lists of elements $\kappa_{\iota}$ $\in$ $\Zc$ such that $\sum_{\iota=1}^{\omega}\kappa_{\iota} \cdot \rho_{\iota}$ $=$ $0$ in $\Zc$, and we sometimes abuse notation and use the last equation interchangeably with the syzygy $(\kappa_1, \dots, \kappa_{\omega})$.  Of course, such syzygies yield $\gf[\mathbf{c}]$-syzygies under the canonical map $\fna{\Zc}{\gf[\mathbf{c}]}$.   

We record some properties of the Jacobi identity syzygies in the following lemma, but first we pause to define the \textbf{spine} of a syzygy $(\kappa_1, \dots, \kappa_{\omega})$ to be 
\begin{equation}\label{eqn:spine}
  \{\rho_{\iota}\, \mid\, \kappa_{\iota} \in \mathbb{Z} \text{ and } \kappa_{\iota} \neq 0\}. 
\end{equation}

\begin{lem} \label{lem:Jsyzprops}
\noindent
  \begin{enumerate}
     \item For the Jacobi identity syzygy $J^{klm}_{pq}$ $=$ $(\kappa_1, \dots, \kappa_{\iota})$, we have that every summand $\kappa_{\iota}\cdot \rho_{\iota}$ in the \emph{LHS} of \emph{(\ref{eqn:Jsyzexp2})} is homogeneous of multi-degree 
\[
    \inc_k(\inc_l(\inc_m(\mdeg(t_q)))) - \mdeg(t_p).
\]
    \item $\sum_{i=1}^{\mu}J^{klm}_{ii}$ $=$ $0$.  
    \item The spine of the syzygy $J^{klm}_{pq}$ comprises $\leq$ $6$ elements, and the components $\kappa_{\iota}$ of $J^{klm}_{pq}$ asociated to the $\rho_{\iota}$ in the spine are all equal to $\pm 1$.  
 \end{enumerate}
 \end{lem}

\emph{Proof:}
The first assertion is a straightforward consequence of Lemma \ref{lem:rhoMultDeg}, since each summand $\kappa_{\iota}\cdot \rho_{\iota}$ is
obtained by combining one or more terms on the LHS of (\ref{eqn:Jsyzexp2}), and the latter is a sum of expressions of the form
\[
     (\pm 1) \cdot e_p \cdot (\text{product of } \gfmm_k, \gfmm_l, \gfmm_m \text{ in some order}) \cdot e_q^{tr}.
\]

To prove the second assertion, view the  $\rho$'s in the matrix 
\[
    {[} \gfmm_k, (\rho^{lm}_{pq})] - [\gfmm_l, (\rho^{km}_{pq})] + [\gfmm_m, (\rho^{kl}_{pq})]
\]
as indeterminates.  Since this matrix is a sum of commutators, its trace is $0$; therefore, the coefficient of every indeterminate $\rho^{k'l'}_{pq}$ in this trace is is $0$.  In particular, this is so for the indeterminates corresponding to the non-trivially-$0$ $\rho_{\iota}$, for which the (vanishing) coefficients in the trace clearly equal the corresponding components of $\sum_{i=1}^{\mu}J^{klm}_{ii}$.

The third assertion is immediate from (\ref{eqn:Jsyzexp2}).  \qedsymbol

\section{Examples of Jacobi identity syzygies} \label{sec:JISyzExs}

We present three examples: the first two show that the syzygies $J^{klm}_{pq}$ simplify drastically in certain cases, and the third will be further studied in section \ref{sec:TSyzEx2}.

\begin{exmp} \label{exmp:JacId1}
    Suppose that $\script{O}$ $\subseteq$ $\mathbb{T}^n$ is such that for some $t_q$ $\in$ $\script{O}$ and choice of variables $x_k$, $x_l$, $x_m$ (with $k$ $<$ $l$ $<$ $m$), we have that  $x_k x_l t_q$ $\in$ $\script{O}$, $x_k x_m t_q$ $\in$ $\script{O}$, and $x_l x_m t_q$ $\in$ $\script{O}$.  Then $J^{klm}_{pq}$ $=$ $0$ for all $1$ $\leq$ $p$ $\leq$ $\mu$.  To see this, we must check that every non-trivially-$0$ $\rho_{\iota}$ on the LHS of (\ref{eqn:Jsyzexp2}) has a coefficient of $0$; by symmetry, it suffices to do this for the $\rho$'s on the first row, which we proceed to examine term-by-term:
\begin{itemize}
  \item ($\rho^{lm}_{\tau_k^{\prime}(p),q}$): note that this $\rho$ is $0$ if $\tau_k^{\prime}(p)$ $=$ $0$, by our conventions (\ref{eqn:conv}), but even if $\tau_k^{\prime}(p)$ $\neq$ $0$, this $\rho$ is trivially $0$ by the second assertion in Corollary \ref{cor:cDegs} and our hypothesis;
  \item ($c_{p,\sigma_k(i)}\rho^{lm}_{iq}$): the $\rho$ in this term is also trivially $0$ by the second assertion in Corollary \ref{cor:cDegs} and our hypothesis;
  \item ($\rho^{lm}_{p,\tau_k(q)}$): this $\rho$ is trivially $0$ again by Corollary \ref{cor:cDegs} and our hypothesis (note that $\tau_k(q)$ $\neq$ $0$);
  \item ($c_{i,\sigma_k(q)}\rho^{lm}_{pi}$): here it is possible that the $\rho$ is not trivially $0$; however, because $\sigma_k(q)$ $=$ $0$ by hypothesis, the coefficient $c_{i,\sigma_k(q)}$ $=$ 0, according to (\ref{eqn:conv}).
\end{itemize}

The simplest nontrivial case of this example is the order ideal 
\[
  \script{O} = \{ 1,\, x_1,\, x_2,\, x_3,\, x_1 x_2,\, x_1 x_3,\, x_2 x_3 \} \subseteq \mathbb{T}^3;
\]
the analysis above shows that $J^{1,2,3}_{p,1}$ $=$ $0$ for $1$ $\leq$ $p$ $\leq$ $7$.
\end{exmp}

\begin{exmp} \label{exmp:JacId2}
More interesting is the case in which $x_k x_l t_q$ $\in$ $\script{O}$ and $x_k x_m t_q$ $\in$ $\script{O}$, but $x_l x_m t_q$ $=$ $b_{j_1}$ $\in$ $\partial \script{O}$.  In this case, one checks that all the terms on the LHS of (\ref{eqn:Jsyzexp2}) are equal to $0$ except for $\rho^{km}_{p,\tau_l(q)}$ and $\rho^{kl}_{p,\tau_m(q)}$, so the syzygy degenerates to 
\[
    \rho^{kl}_{p,\tau_m(q)}\ =\ \rho^{km}_{p,\tau_l(q)}.
\]
This relation can be derived directly by observing that $\rho^{kl}_{p,\tau_m(q)}$ and $\rho^{km}_{p,\tau_l(q)}$ are ``Case 3'' $\rho$'s (see Section \ref{sec:BBSidGens}) involving the same pair of monomials in $\partial \script{O}$: 
\[
  b_{j_1} = x_l \cdot t_{\tau_m(q)} = x_m \cdot t_{\tau_l(q)},\ \text{and } b_{j_2} = x_k \cdot b_{j_1},
\]
so the equality follows from equation (\ref{eqn:Case3}).

The simplest case of this example is 
\[
  \begin{array}{c}
    \script{O} = \{ 1,\, x_1,\, x_2,\, x_3,\, x_1 x_2,\, x_1 x_3 \} \subseteq \mathbb{T}^3,\vspace{.05in}\\
    \partial \script{O} = \{ x_1^2,\, x_2^2,\, x_2 x_3,\, x_3^2,\, x_1^2 x_2,\, x_1 x_2^2,\, x_1 x_2 x_3,\, x_1^2 x_3,\, x_1 x_3^2 \}.
  \end{array}
\]
We index the $t_i$ $\in$ $\script{O}$ and the $b_j$ $\in$ $\partial \script{O}$ as they are displayed.  Then for $t_1$ $=$ $1$, one checks easily that $x_1 x_2 t_1$ $\in$ $\script{O}$ and $x_1 x_3 t_1$ $\in$ $\script{O}$, but $x_2 x_3 t_1$ $\notin$ $\script{O}$.  Therefore, we have that $\rho^{1,3}_{p,3}$ $=$ $\rho^{1,2}_{p,4}$ for $1$ $\leq$ $p$ $\leq$ $6$; for example, one checks that 
\[
  \begin{array}{c}
    \rho^{1,3}_{1,3}\ =\ \rho^{1,2}_{1,4}\ =\ -c_{1,7} + c_{1,1}\, c_{2,3} + c_{1,5}\,
   c_{5,3} + c_{1,8}\, c_{6,3},\vspace{.05in}\\
    \rho^{1,3}_{6,3}\ =\ \rho^{1,2}_{6,4}\ =\ c_{4,3} + c_{2,3}\, c_{6,1} + c_{5,3}\,
   c_{6,5} - c_{6,7} + c_{6,3}\, c_{6,8}.
  \end{array}
\]
\end{exmp}

\begin{exmp} \label{exmp:JacId3}
    Let $\script{O}$ $=$ $\{1, x_1\}$ $\subseteq$ $\mathbb{T}^3$.  One sees easily that 
\[
    \partial \script{O} = \{ x_2,\ x_3,\ x_1^2,\ x_1 x_2,\ x_1 x_3\},
\]
and that the three multiplication matrices are
\begin{equation} \label{eqn:JExMultMats}
  \gfmm_1 = \left( \begin{array}{cc}
                    0 & c_{1,3}\\
                    1 & c_{2,3}
                   \end{array} \right),\
 \gfmm_2 = \left( \begin{array}{cc}
                    c_{1,1} & c_{1,4}\\
                    c_{2,1} & c_{2,4}
                   \end{array} \right),\ 
 \gfmm_3 = \left( \begin{array}{cc}
                    c_{1,2} & c_{1,5}\\
                    c_{2,2} & c_{2,5}
                   \end{array} \right).\
\end{equation}
The commutators of the multiplication matrices are then
\begin{equation} \label{eqn:JExComms}
    \begin{array}{c}
      [\gfmm_1, \gfmm_2]\ =\ \left(\begin{array}{cc}
                                 \rho^{1,2}_{1,1} & \rho^{1,2}_{1,2}\vspace{.05in}\\
                                 \rho^{1,2}_{2,1} & \rho^{1,2}_{2,2}
                               \end{array}\right)\ = \vspace{.05in}\\
\left(
  \begin{array}{cc}
 c_{1,3}\, c_{2,1}-c_{1,4} & -c_{1,1}\,
   c_{1,3}+c_{2,4}\, c_{1,3}-c_{1,4}\, c_{2,3}\vspace{.05in} \\
 c_{1,1}+c_{2,1}\, c_{2,3}-c_{2,4} &
   c_{1,4}-c_{1,3}\, c_{2,1}
  \end{array}
\right),\vspace{.1in}\\

[\gfmm_1, \gfmm_3]\ =\ \left(\begin{array}{cc}
                                 \rho^{1,3}_{1,1} & \rho^{1,3}_{1,2}\vspace{.05in}\\
                                 \rho^{1,3}_{2,1} & \rho^{1,3}_{2,2}
                               \end{array}\right)\ = \vspace{.05in}\\
\left(
\begin{array}{cc}
 c_{1,3}\, c_{2,2}-c_{1,5} & -c_{1,2}\,
   c_{1,3}+c_{2,5}\, c_{1,3}-c_{1,5}\, c_{2,3}\vspace{.05in} \\
 c_{1,2}+c_{2,2}\, c_{2,3}-c_{2,5} &
   c_{1,5}-c_{1,3}\, c_{2,2}
\end{array}
\right), \vspace{.1in}\\

[\gfmm_2, \gfmm_3]\  =\ \left(\begin{array}{cc}
                                 \rho^{2,3}_{1,1} & \rho^{2,3}_{1,2}\vspace{.05in}\\
                                 \rho^{2,3}_{2,1} & \rho^{2,3}_{2,2}
                               \end{array}\right)\ = \vspace{.05in}\\
\left(
\begin{array}{cc}
 c_{1,4}\, c_{2,2}-c_{1,5}\, c_{2,1} & -c_{1,2}\,
   c_{1,4}+c_{2,5}\, c_{1,4}+c_{1,1}\,
   c_{1,5}-c_{1,5}\, c_{2,4} \\
 c_{1,2}\, c_{2,1}-c_{2,5}\, c_{2,1}-c_{1,1}\,
   c_{2,2}+c_{2,2}\, c_{2,4} & c_{1,5}\,
   c_{2,1}-c_{1,4}\, c_{2,2}
\end{array}
\right);
    \end{array}
\end{equation}
note that all of the $\rho$'s are non-zero in this case.

The Jacobi identity syzygies are given by the four entries of the matrix
\[
    [\gfmm_1, [\gfmm_2, \gfmm_3]] - [\gfmm_2, [\gfmm_1, \gfmm_3]] + [\gfmm_3, [\gfmm_1, \gfmm_2]].
\]
Before listing the matrix entries, we recall that the $(1, 1)$-entry and the $(2, 2)$-entry are negations of one another, by the second assertion in Lemma \ref{lem:Jsyzprops}, so we only have three distinct syzygies:
\begin{equation} \label{eqn:JsyzEg3list}
  \begin{array}{rcl}
    J^{1,2,3}_{1,1} & = & \left( -c_{2,2}\, \rho ^{1,2}_{1,2}+c_{1,5}\, \rho
   ^{1,2}_{2,1} + c_{2,1}\, \rho^{1,3}_{1,2}-c_{1,4}\, \rho
  ^{1,3}_{2,1}+c_{1,3}\, \rho^{2,3}_{2,1} - \rho^{2,3}_{1,2}\ =\ 0\right),\vspace{.1in}\\
    J^{1,2,3}_{1,2} & = & \left( \left\{\begin{array}{l}
       -c_{1,5}\, \rho^{1,2}_{1,1} + (c_{1,2} - c_{2,5})\, \rho^{1,2}_{1,2}       + c_{1,5}\, \rho
       ^{1,2}_{2,2} + c_{1,4}\, \rho^{1,3}_{1,1}\ + \vspace{.05in}\\
       (c_{2,4} - c_{1,1})\, \rho
       ^{1,3}_{1,2}\ - c_{1,4}\, \rho
       ^{1,3}_{2,2} - c_{1,3}\, \rho^{2,3}_{1,1} - c_{2,3}\, \rho
       ^{2,3}_{1,2} + c_{1,3}\, \rho^{2,3}_{2,2}
      \end{array}\right\} = 0 \right),\vspace{.1in}\\
   J^{1,2,3}_{2,1} & = & \left( \left\{\begin{array}{l} c_{2,2}\, \rho^{1,2}_{1,1} +(c_{2,5} - c_{1,2})\, \rho
   ^{1,2}_{2,1} - c_{2,2}\, \rho
   ^{1,2}_{2,2} - c_{2,1}\, \rho^{1,3}_{1,1}\ + \vspace{.05in}\\
   (c_{1,1} - c_{2,4})\, \rho
   ^{1,3}_{2,1} + c_{2,1}\, \rho
   ^{1,3}_{2,2} + c_{2,3}\, \rho^{2,3}_{2,1} + \rho
   ^{2,3}_{1,1} - \rho^{2,3}_{2,2} \end{array}\right\}  = 0 \right).
  \end{array}
\end{equation}

These syzygies can be further simplified by using the relations
\[
    \rho^{1,2}_{2,2} = -\rho^{1,2}_{1,1},\ \ \rho^{1,3}_{2,2} = -\rho^{1,3}_{1,1},\ \ \rho^{2,3}_{2,2} = -\rho^{2,3}_{1,1},
\]
which are special cases of the ``trace syzygies'' to be defined in section \ref{sec:Tsyz}.  The trace syzygies are based on the fact that the trace of a commutator of two square matrices is $0$.  In the next section, we lay the groundwork for their definition by studying certain expressions that reduce to commutators.
\end{exmp}

\section{Sums of products that reduce to commutators} \label{sec:sumOfCom}

	Let $m_1, \dots, m_n$ (for $n$ $\geq$ $2$) be non-commuting indeterminates.  In this section we exhibit certain sums of products of the $m_k$ and the basic commutators $[m_k, m_l]$ $=$ $(m_k m_l - m_l m_k)$ that simplify to single commutators.  In the sequel, the $m_k$ will be replaced by the generic formal multiplication matrices $\gfmm_k$ associated to a given order ideal $\script{O}$; by taking the trace of each of our expressions, we will obtain syzygies of the $\rho^{kl}_{pq}$, which are the entries of the matrices $[\gfmm_k, \gfmm_l]$, $1$ $\leq$ $k$ $<$ $l$ $\leq$ $n$.

	To produce an expression of the desired type, we begin by choosing 
an ordered product
\begin{equation} \label{eqn:ordProd}
  \Pi\ =\ m_{k_1} m_{k_2} \dots m_{k_{s}}
\end{equation}
such that at least two of the indeterminates (say $m_k$ and $m_l$) appear at least once.  We define the multi-degree of each $m_k$ to be $\mdeg(m_k)$ $=$ $e_k^{tr}$, and extend the definition to $\Pi$ in the obvious way, so that $\mdeg(\Pi)$ $=$ $(d_1, \dots, d_n)^{tr}$ is a column of non-negative integers such that at least two of the components $d_k$ and $d_l$ are positive; such ordered products and tuples will be called \textbf{good}, or, more precisely, $(k,l)$\textbf{-good}.  Next, we choose one of the indices $k$ for which the associated integer $d_k$ $>$ $0$, and call this the \textbf{distinguished index}; the variable $m_k$ is then called the \textbf{distinguished variable}.  We then form a shorter ordered product $\Pi_{\hat{k}}$ by deleting one of the occurrences of $m_k$ from the product $\Pi$ (it doesn't matter which occurrence of $m_k$ is removed, but for definiteness we will remove the left-most).  For example, if $n$ $=$ $3$, 
\[
    \Pi\ =\ m_2 m_1 m_3 m_1 m_2\ \Rightarrow\ \mdeg(\Pi)\ =\ (d_1, d_2, d_3)^{tr}\ =\ (2, 2, 1)^{tr}, 
\]
and $1$ is the distinguished index, then
\[
    \Pi_{\hat{1}}\ =\ m_2 m_3 m_1 m_2.
\]

We now  sum the terms obtained by replacing, one at a time, each of the $m_l$ in the ordered product $\Pi_{\hat{k}}$ by the basic commutator $[m_k,m_l]$.  For example, the expression formed in this way from $\Pi_{\hat{1}}$ above is 
\[
    [m_1, m_2] m_3 m_1 m_2 + m_2 [m_1, m_3] m_1 m_2  + m_2 m_3 [m_1, m_1] m_2 + m_2 m_3 m_1 [m_1, m_2].
\]
(Our assumption that at least two of the variables $m_k$ and $m_l$ actually appear in $\Pi$ ensures that we avoid the trivial case in which every commutator in our expression is equal to $[m_k, m_k]$ $=$ $0$.)
It is easily verified by hand that the preceding expression simplifies to a commutator (namely, $[m_1, m_2 m_3 m_1 m_2]$); the following proposition generalizes this observation.

\begin{prop} \label{prop:mainCommLem}
	Let $\Pi$ be a $(k,l)$-good ordered product with $\mdeg(\Pi)$ $=$ $(d_1, \dots, d_n)^{tr}$, and $\Pi_{\hat{k}}$ $=$ $m_{l_1} m_{l_2} \dots m_{l_{s-1}}$ the shorter ordered product obtained from $\Pi$ as described above, with $k$ the distinguished index.  Then the expression
\[
	\sum_{v=1}^{s-1} \left(m_{l_1} m_{l_2} \dots m_{l_{v-1}}[m_k, m_{l_{v}}] m_{l_{v+1}} \dots m_{l_{s-1}}\right)
\]
is homogeneous of multi-degree $(d_1,\dots, d_n)^{tr}$ in the variables $m_l$, $1$ $\leq$  $l$ $\leq$ $n$, and
reduces to the commutator $[ m_k,\ m_{l_1} m_{l_2} \dots m_{l_{s-1}}]$.
\end{prop}

\emph{Proof:} The first statement is obvious.  For the second, it suffices to observe that the negative contribution from the $v$-th term of the sum:
\[
	 -\ m_{l_1} m_{l_2}\dots m_{l_{v-1}} (m_{l_{v}}\cdot m_k) m_{l_{v+1}} \dots m_{l_{s-1}},
\]
cancels with the positive contribution from the $v+1$-st term:
\[
	+\ m_{l_1}m_{l_2} \dots m_{l_{v}}(m_k \cdot m_{l_{v}+1}) m_{k_{v+2}}\dots m_{l_{s-1}},
\]
so the sum telescopes to yield
\[
  \left( (m_k \cdot m_{l_1})m_{l_2} \dots m_{l_{s-1}}\right) 
- 
 \left( m_{l_1} \dots m_{l_{s-2}}(m_{l_{s-1}} \cdot m_k) \right)\ =\ [m_k, m_{l_1}\dots m_{l_{s-1}}], 
\]
as desired. \qedsymbol
\bigskip

It is convenient to adopt the following ``index-tuple notation'' for ordered products:
\begin{equation} \label{eqn:indNotn}
    \Pi\ =\ m_{k_1} m_{k_2} \dots m_{k_{s}}\ =\  \left< k_1, k_2, \dots, k_s \right> .
\end{equation}

\section{The trace syzygies} \label{sec:Tsyz}

     The basic idea is very simple.  We first replace the indeterminates $m_l$ used in the preceding section by the generic multiplication matrices $\gfmm_l$ associated to the order ideal $\script{O}$ $=$ $\{t_1, \dots, t_{\mu} \}$.  Then, given a $(k,l)$-good ordered product 
\[
    \Pi\ =\ \gfmm_{k_1}\, \gfmm_{k_2}\, \dots,\, \gfmm_{k_s}\ =\ \left< k_1, k_2, \dots, k_s\right>
\]
with multi-degree $\mdeg(\Pi)$ $=$ $(d_1, \dots, d_n)^{tr}$, and choice of distinguished index $k$, we extract the shorter product $\Pi_{\hat{k}}$, and observe that the expression formed as in Proposition \ref{prop:mainCommLem}:
\[
	\sum_{v=1}^{s-1} \gfmm_{l_1}\, \gfmm_{l_2}\dots \gfmm_{l_{v-1}}\, [\gfmm_k, \gfmm_{l_{v}}]\, \gfmm_{l_{v+1}}\dots \, \gfmm_{l_{s-1}},
\]
since it reduces to a commutator of two matrices, has trace 0.  Furthermore, the entries of the matrices $\gfmm_l$ are elements of the set 
\[
  \{0, 1\}\ \cup\ \{ c_{ij} \mid 1 \leq i \leq \mu,\ 1 \leq j \leq \nu \}, 
\]
and the entries of the commutators $[\gfmm_k, \gfmm_{l_{v}}]$ are either equal to $\rho^{kl_{v}}_{pq}$ if $k$ $<$ $l_v$, $0$ if $k$ $=$ $l_v$, or $-\rho^{l_{v}k}_{pq}$ if $k$ $>$ $l_v$.  Hence, the trace is a $\Zc$-linear combination of the $\rho$'s that is equal to 0, and so yields a syzygy thereof.  As in Section \ref{sec:JISyz}, we ignore the $\rho$'s that are trivially $0$, and use the same list $\rho_1, \rho_2, \dots, \rho_{\omega}$ of non-trivially-$0$ $\rho$'s.  Accordingly, we define the \textbf{trace syzygy} having \textbf{ordered product} 
$\Pi$ $=$  $\left< k_1, k_2, \dots k_s\right>$ and \textbf{distinguished index} $k$ to be the sequence of coefficients 
\[
   T\ =\ T_{\Pi,k}\ =\ (\kappa_1, \dots, \kappa_{\omega}), 
\] 
where $\kappa_{\iota}$ $\in$ $\Zc$ is the coefficient of $\rho_{\iota}$ in the expression
\begin{equation} \label{eqn:TSyzSum}
    \Tr \left(		\sum_{v=1}^{s-1} \gfmm_{l_1}\, \gfmm_{l_2}\dots \gfmm_{l_{v-1}}\, [\gfmm_k, \gfmm_{l_{v}}]\, \gfmm_{l_{v+1}}\dots \, \gfmm_{l_{s-1}} \right).
\end{equation}
Again, we will abuse the notation and write this syzygy as $\sum_{\iota=1}^{\omega} \kappa_{\iota} \cdot \rho_{\iota}$ $=$ $0$.

\begin{lem} \label{lem:TSyzMulDeg}
    The expression \emph{(\ref{eqn:TSyzSum})} is homogeneous of multi-degree
\[
    \mdeg(\Pi)\ =\  (d_1, \dots, d_n)^{tr};
\]
whence, each polynomial $\kappa_{\iota} \cdot \rho_{\iota}$ has the same property. 
\end{lem}

\emph{Proof:}
    The expression (\ref{eqn:TSyzSum}) is a sum of contributions coming from the traces of products of the form 
\[
    \gfmm_{l_1}\, \gfmm_{l_2}\dots \gfmm_{l_{v-1}}\, [\gfmm_k, \gfmm_{l_{v}}]\, \gfmm_{l_{v+1}}\dots \, \gfmm_{l_{s-1}}.  
\]
By Lemma \ref{lem:rhoMultDeg}, we know that if the $q$-th elementary unit column vector $e_q^{tr}$ (of length $\mu$) is multiplied by a product of the matrices $\gfmm_{l_1}$, \dots, $\gfmm_{l_{s-1}},$ and $\gfmm_k$ in some order, the result is a column vector $(\beta_1, \dots, \beta_{\mu})$ such that the polynomial $\sum_{i=1}^{\mu}\beta_i t_i$ is homogeneous of multi-degree $\mdeg(t_q)$ $+$ $(d_1, \dots, d_n)^{tr}$.  Therefore, $\beta_q$, which is the $(q,q)$-th entry of the matrix product, is homogeneous of multi-degree
\[
    \mdeg(t_q) + (d_1, \dots, d_m)^{tr} - \mdeg(t_q)\ =\ (d_1, \dots, d_m)^{tr};
\]
the desired result follows immediately from this.\  \qedsymbol
\medskip

We proceed to study the spine (\ref{eqn:spine}) of a trace syzygy $T_{\Pi, k}$ $=$ $(\kappa_1, \dots, \kappa_{\omega})$.  The next result is the heart of the matter:

\begin{prop} \label{prop:Tspines}
 
     Let $T$ $=$ $T_{\Pi,k}$ be the trace syzygy having ordered product $\Pi$ $=$ $\left<k_1, k_2, \dots, k_s \right>$ and distinguished index $k$, with $\mdeg(\Pi)$ $=$ $(d_1, \dots, d_n)^{tr}$ a $(k,l)$-good tuple for at least one index $l$ $\neq$ $k$).  Choose one of the non-trivially-$0$ $\rho$'s, say $\rho_{\iota}$ $=$ $\rho^{k'l'}_{pq}$.  Then the following are equivalent:
\begin{enumerate}
  \item The $\iota$-th component $\kappa_{\iota}$ of $T$ is a non-zero constant $\kappa_{\iota}$ $\in$ $\mathbb{Z}$; that is, $\rho_{\iota}$ is in the spine of $T$.
  \item The $\iota$-th component $\kappa_{\iota}$ of $T$, which is an element of $\Zc$, has a non-zero constant term $\kappa_{\iota,0}$ $\in$ $\mathbb{Z}$.
  \item $\{k',l'\}$ $=$ $\{ k,l \}$ for some $l$ $\in$ $\{ k_1, \dots, k_s \} \setminus \{ k \}$, and $\rho^{k'l'}_{pq}$ is homogeneous of multi-degree
\[
    \mdeg(\rho^{k'l'}_{pq})\ =\ \mdeg(\Pi)\ =\ (d_1, \dots, d_n)^{tr}.
\]
\end{enumerate}
Furthermore, if any of the equivalent conditions holds, we have that $\kappa_{\iota}$ is a non-zero element of $\mathbb{Z}$ given by
\[
    \kappa_{\iota} = \kappa_{\iota,0} =
      \left\{  \begin{array}{l}
            d_{l'} \text{ if } k = k'\\
           -d_{k'} \text{ if } k = l'.
              \end{array}\right. 
\]
\end{prop}

\emph{Proof:}
(1) $\Rightarrow$ (2): Immediate.

(2) $\Rightarrow$ (3): Recall that $\kappa_{\iota}$ is the coefficient of $\rho_{\iota}$ $=$ $\rho^{k'l'}_{pq}$ in the expression (\ref{eqn:TSyzSum}).  We are assuming that $\kappa_{\iota}$ has a non-zero constant term, so it must be the case that for some index $1$ $\leq$ $v$ $\leq$ $s-1$, the coefficient of $\rho^{k'l'}_{pq}$ in the $v$-th summand
\begin{equation}\label{eqn:vthSummand}
    \Tr \left( \gfmm_{l_1}\, \gfmm_{l_2}\dots \gfmm_{l_{v-1}}\, [\gfmm_k, \gfmm_{l_{v}}]\, \gfmm_{l_{v+1}}\dots \, \gfmm_{l_{s-1}} \right)
\end{equation}
of (\ref{eqn:TSyzSum}) has a non-zero constant term $\lambda_0$.  Thus $\rho^{k'l'}_{pq}$ must appear, either positively or negatively, as one of the entries of the commutator $[\gfmm_k, \gfmm_{l_v}]$, whose entries are of the form $\rho^{kl_v}_{pq}$ if $k$ $<$ $l_v$ and $-\rho^{l_vk}_{pq}$ if $k$ $>$ $l_v$ (clearly $k$ $=$ $l_v$ is impossible here).  It follows immediately that $\{k', l'\}$ $=$ $\{k, l_v\}$, with $l_v$ $=$ $l$ $\in$ $\{ k_1, \dots, k_s \} \setminus \{ k \}$.  Furthermore, since $\lambda_{0}\cdot \rho^{k'l'}_{pq}$ is part of the expression (\ref{eqn:TSyzSum}), we have by Lemma \ref{lem:TSyzMulDeg} that $\mdeg(\lambda_{0}\cdot \rho^{k'l'}_{pq})$ $=$ $\mdeg(\rho^{k'l'}_{pq})$ $=$ $(d_1, \dots, d_n)^{tr}$.

(3) $\Rightarrow$ (1): Suppose that we are in case $k'$ $=$ $k$, so that $l'$ $=$ $l$ $=$ $l_v$ in the indexing of (\ref{eqn:TSyzSum}) --- there must exist at least one such index.  \emph{We claim} that $\rho^{kl_v}_{pq}$ appears in the trace of the $v$-th summand (\ref{eqn:vthSummand}) with coefficient $1$.  (It will be clear that the coefficient is $-1$ in the other case, when $l'$ $=$ $k$.)  

To prove the claim, first note that the following equations
\[
    \mdeg(\rho^{kl_v}_{pq})\ =\ \inc_k(\inc_{l_v}(\mdeg(t_q))) - \mdeg(t_p)\ =\  (d_1, \dots, d_n)^{tr}\ =\ \mdeg(\Pi),
\]
which come from Lemma \ref{lem:rhoMultDeg} and our hypothesis,
imply that
\[
    \frac{x_k\, x_{l_v}\, t_q}{x_{k_1}\, x_{k_2}\,\cdots\, x_{k_s}} \ =\ \frac{t_q}{x_{l_1}\, x_{l_2}\,\cdots\, x_{l_{v-1}}\, x_{l_{v+1}}\, \cdots\, x_{l_{s-1}}}\ =\ t_p.
\]
Let $1$ $\leq$ $r$ $\leq$ $\mu$ be the index such that 
\begin{equation} \label{eqn:trDef}
  t_r = (x_{l_1} x_{l_2}\cdots x_{l_{v-1}})\cdot t_p\ \ \Rightarrow\ \  (x_{l_{v+1}} x_{l_{v+2}}\cdots x_{l_{s-1}})\cdot t_r = t_q.
\end{equation}
Any term involving $\rho^{kl_v}_{pq}$ in (\ref{eqn:vthSummand})  would have to appear in the $(u,u)$-th component of the matrix product for some $1$ $\leq$ $u$ $\leq$ $\mu$, which is given by
\[
    e_u \cdot \left( \gfmm_{l_1}\, \gfmm_{l_2}\dots \gfmm_{l_{v-1}}\, [\gfmm_k,  \gfmm_{l_v}]\, \gfmm_{l_{v+1}}\dots \, \gfmm_{l_{s-1}} \right) \cdot e_u^{tr},
\]
where $e_u$ denotes the $u$-th standard unit vector.  Suppose that
\[
    \left( \gfmm_{l_{v+1}}\dots \, \gfmm_{l_{s-1}} \right) \cdot e_u^{tr}\ =\ \mathbf{\gamma}^{tr} = (\gamma_1, \gamma_2, \dots, \gamma_{\mu})^{tr},\ \ \gamma_i \in \Zc. 
\]
Then $\rho^{kl_v}_{pq}$ can appear in the $p$-th component (only) of
\[
    [\gfmm_k, \gfmm_{l_v})] \cdot \mathbf{\gamma}^{tr}\ =\ (\rho^{kl_v}_{ij})\cdot \mathbf{\gamma}^{tr},
\]
which is equal to $\sum_{j=1}^{\mu}\rho^{kl_v}_{pj}\cdot \gamma_j$.  It remains to compute 
\[
  \begin{array}{cl}
   {} &  e_u \cdot \left( \gfmm_{l_1}\, \gfmm_{l_2}\dots \gfmm_{l_{v-1}} \right) \cdot (\sum_{j=1}^{\mu}\rho^{kl_v}_{pj}\cdot \gamma_j)e_p^{tr}\vspace{.05in}\\
 = & (\sum_{j=1}^{\mu}\rho^{kl_v}_{pj}\cdot \gamma_j) \cdot e_u \cdot \left( \gfmm_{l_1}\, \gfmm_{l_2}\dots \gfmm_{l_{v-1}} \right) \cdot e_p^{tr}.
  \end{array}
\]
But, as noted above,  $(x_{l_1} x_{l_2} \cdots x_{l_{v-1}})\cdot t_p$ $=$ $t_r$ $\in$ $\script{O}$, so the above product reduces to 
\[
    (\sum_{j=1}^{\mu}\rho^{kl_v}_{pj}\cdot \gamma_j) \cdot e_u  \cdot e_r^{tr},
\]
which is equal to $0$ unless $u$ $=$ $r$.  In other words, $\rho^{kl_v}_{pq}$ can only appear in the trace by appearing in the $(r,r)$-th component of the matrix product, which equals 
\[
  \begin{array}{cl}
     {} & e_r \cdot \left( \gfmm_{l_1}\, \gfmm_{l_2}\dots \gfmm_{l_{v-1}}\, [\gfmm_k, \gfmm_{l_v}]\, \gfmm_{l_{v+1}}\dots \, \gfmm_{l_{s-1}} \right) \cdot e_r^{tr}\vspace{.05in}\\
     = &   e_r \cdot \left( \gfmm_{l_1}\, \gfmm_{l_2}\dots \gfmm_{l_{v-1}}\right) \cdot (\rho^{kl_v}_{ij}) \cdot e_q^{tr}\ \  \text{(by (\ref{eqn:trDef}))} \vspace{.05in}\\
     = & e_r \cdot \left( \gfmm_{l_1}\, \gfmm_{l_2}\dots \gfmm_{l_{v-1}}\right) \cdot (\rho^{kl_v}_{1,q},\, \rho^{kl_v}_{2,q},\, \dots,\, \rho^{kl_v}_{\mu,q})^{tr},
  \end{array}
\]
and the only term involving $\rho^{kl_v}_{pq}$ in this expression is
\[
  \begin{array}{rcl}
    e_r \cdot \left( \gfmm_{l_1}\, \gfmm_{l_2}\dots \gfmm_{l_{v-1}}\right) \cdot \rho^{kl_v}_{pq} e_p^{tr} &  = & \rho^{kl_v}_{pq} \left( e_r \cdot \left( \gfmm_{l_1}\, \gfmm_{l_2}\dots \gfmm_{l_{v-1}}\right) \cdot e_p^{tr}\right)\vspace{.05in}\\ {} & = &  \rho^{kl_v}_{pq} (e_r \cdot e_r^{tr})\vspace{.05in}\\
{} & = & \rho^{kl_v}_{pq}.
  \end{array}
\]
This completes the proof of the claim.

To finish the proof of (1), we compute $\kappa_{\iota}$ by summing the coefficients of $\rho^{k'l'}_{pq}$ in each of the summands (\ref{eqn:vthSummand}), $1$ $\leq$ $v$ $\leq$ $s-1$.  Since this coefficient can only be non-zero in case $\{k', l'\}$ $=$ $\{k, l_v \}$, and for any such $v$, of which there are $d_{l_v}$ in all, this coefficient is $+1$ if $k$ $=$ $k'$, and $-1$ if $k$ $=$ $l'$, we obtain 
\begin{equation}\label{eqn:kappaiota}
    \kappa_{\iota}\ =\ \left\{ \begin{array}{l}
                          d_{l_v}  =  d_{l'}, \text{ if } k = k'\\
                          -d_{l_v} =   -d_{k'}, \text{ if } k = l'
                       \end{array}; \right.
\end{equation}
in particular, we have that $\kappa_{\iota}$ $\in$ $\mathbb{Z}$, as desired. This completes the equivalence proof.   

Finally, it is now clear that each of our equivalent conditions entails that $\kappa_{\iota}$ is a non-zero integer whose value is given by (\ref{eqn:kappaiota}); whence, the last assertion in the proposition holds, and we are done.
\qedsymbol
 
\begin{cor} \label{cor:T1}
    If the component $\kappa_{\iota}$ of $T_{\Pi,k}$ is not a constant, then its constant term is $0$.  \qedsymbol
\end{cor}

\emph{Proof:}
    The proposition shows that if $\kappa_{\iota}$ has a non-zero constant term, then it is in fact a constant.
\qedsymbol

\begin{cor} \label{cor:T3}
    Let $\Pi$ $=$ $\left<k_1, k_2, \dots, k_s\right>$ be a good ordered product of multi-degree $\mdeg(\Pi)$ $=$ $(d_1, \dots, d_n)^{tr}$, and let $S$ $=$ $\{ k\, \mid\, d_k > 0  \}$.  Then the following linear combination has empty spine; that is, it has no non-zero constant entries: 
\[
    \sum_{k \in S}  d_k \cdot T_{\Pi,k}.
\]
\end{cor}

\emph{Proof:} Corollary \ref{cor:T1} shows that the entries of $T_{\Pi, k}$ are either constants or have $0$ as constant term; furthermore, it is clear that the $\iota$-th component of $T_{\Pi,k}$ can only be non-zero in case $\rho_{\iota}$ $=$ $\rho^{k'l'}_{pq}$ with $\{k', l'\}$ $=$ $\{k, l\}$ $\subseteq$ $S$.  Consider two indices $k$, $l$ $\in$ $S$ with $k$ $<$ $l$.  If $\rho^{kl}_{pq}$ $=$ $\rho_{\iota}$ is in the spine of $T_{\Pi, k}$, the proposition shows that ${\iota}$-th component of  $T_{\Pi, k}$ is $d_l$.  Therefore, the $\iota$-th component in $d_k \cdot T_{\Pi, k}$ is $d_k d_l$.  The proposition further yields that $\rho^{kl}_{pq}$ is also in the spine of $T_{\Pi, l}$, and the $\iota$-th component of $T_{\Pi, l}$ is -$d_k$, so the $\iota$-th component of $d_l \cdot T_{\Pi, l}$ is $-d_l d_k$.  These values cancel, and none of the other terms in the linear combination can have a non-zero $\iota$-th component.  Applying this reasoning to all index pairs $k$ $<$ $l$ in $S$, we obtain the desired result.  
\qedsymbol

\begin{cor} \label{cor:T2}
    If $\Pi'$ is a rearrangement of the ordered product $\Pi$, then the tuple $T_{\Pi, k}$ $-$ $T_{\Pi', k}$ has empty spine; in particular, the spines of $T_{\Pi,k}$ and $T_{\Pi',k}$ are equal.
\end{cor}

\emph{Proof:} Since the truth of condition (3) of Proposition \ref{prop:Tspines} depends only on the multi-degree of the ordered product, and not on the order of its factors, we have that $\rho^{k'l'}_{pq}$ is in the spine of $T_{\Pi,k}$ if and only if it is in the spine of $T_{\Pi',k}$; moreover, the corresponding coefficients $\kappa_{\iota}$ and $\kappa'_{\iota}$ are equal.   \qedsymbol
\bigskip

We now describe a procedure for identifying the ordered products $\Pi$ such that $T_{\Pi,k}$ has non-empty spine for some distinguished index $k$.  We will call such ordered products, and their multi-degrees $\mdeg(\Pi)$, \textbf{spinal}.

\newcommand{\arr}{\mathfrak{a}}
We can visualize spinal multi-degrees and the spines of their associated trace syzygies using arrows in the lattice of monomials in $x_1$, $\dots$, $x_n$.  We begin by mapping each $\rho^{kl}_{pq}$ to the arrow $\arr(\rho^{kl}_{pq})$ with tail at $t_p$ and head at $x_k x_l t_q$.  Lemma \ref{lem:rhoMultDeg} yields that
\[
  \mdeg(\rho^{kl}_{pq}) = \inc_k(\inc_l(\mdeg(t_q))) - \mdeg(t_p)
\]
is the displacement (column) vector from the tail to the head of the arrow $\arr(\rho^{kl}_{pq})$.

In order to recognize the non-trivially-$0$ $\rho$'s in terms of their associated arrows, we define a \textbf{target monomial} to be a monomial $m$ for which there exists a pair of indices $1 \leq k < l \leq n$ such that
\[
    m\ =\  x_k\, x_l\, t_q,\ \text{where  } t_q \in \script{O}\ \text{ and  either  }  x_k t_q \notin \script{O} \text{ or } x_l t_q \notin \script{O};
\]
in particular, $m$ $\notin$ $\script{O}$ must hold ($m$ could be, but need not be, an element of $\partial \script{O}$).  More precisely, we call $m$ a $(k,l)$-target monomial if we need to specify the pair of indices that witness to $m$'s target-monomial status; it is of course possible that $m$ could have more than one pair of witnesses, and therefore the mapping $\rho^{kl}_{pq}$ $\mapsto$ $\arr(\rho^{kl}_{pq})$ need not be one-to-one.
One then checks easily that
$\rho^{kl}_{pq}$ is not trivially $0$ if and only if the head of $\arr(\rho^{kl}_{pq})$ is a $(k,l)$-target monomial, since the latter condition is equivalent to $\rho^{kl}_{pq}$ being in Case 3 or Case 4 as described in Section \ref{sec:BBSidGens}.  

Most of Proposition \ref{prop:Tspines} can now be re-expressed in the following, more geometric, fashion:

\begin{cor} \label{cor:TarMD}
Let $\Pi$, $\mdeg(\Pi)$ $=$ $(d_1, \dots, d_n)^{tr}$, etc., be as in the Proposition  \emph{\ref{prop:Tspines}}.  Then $\rho^{k'l'}_{pq}$ is in the spine of $T_{\Pi,k}$ if and only if the following conditions hold:
\begin{itemize}
  \item $\{k', l'  \}$ $=$ $\{ k, l \}$ for some $l$ $\neq$ $k$ with $d_l$ $>$ $0$;
  \item  the arrow $\arr(\rho^{k'l'}_{pq})$ has displacement vector $(d_1, \dots, d_n)^{tr}$;
  \item the head $x_k x_l t_q$ of $\arr(\rho^{k'l'}_{pq})$ is a $(k,l)$-target monomial.
\end{itemize}
  \qedsymbol
\end{cor}

\begin{rem} \label{rem:SpAr}
    In light of the preceding corollary, we have the following simple but tedious-to-describe procedure to determine if a given (good) ordered product $\Pi$, or multi-degree $\mathbf{d}$ $=$ $\mdeg(\Pi)$ $=$ $(d_1, \dots, d_n)^{tr}$, is spinal, and if so, to compute the spines of the trace syzygies $T_{\Pi, k}$, which by definition will be nonempty for at least one $k$.  First, one finds all the arrows $\mathbf{a}$ $=$ $(t_p$ $\rightarrow$ $m)$ such that $t_p$ $\in$ $\script{O}$, $m$ is a target monomial, and the displacement vector is $\mathbf{d}$.  For each such $\mathbf{a}$, one lists the triples 
\[
    (t_{q_{\mathbf{a},\alpha}},\, k_{\mathbf{a},\alpha},\, l_{\mathbf{a},\alpha}),\ \ t_{_{\mathbf{a},\alpha}} \in \script{O},\ \ 0 \leq k_{\mathbf{a},\alpha} < l_{\mathbf{a},\alpha} \leq n,
\]
such that $d_{k_{\mathbf{a},\alpha}} > 0$,  $d_{l_{\mathbf{a},\alpha}} > 0$, and the pair $(k_{\mathbf{a},\alpha},$ $l_{\mathbf{a},\alpha})$ witnesses that $m$ $=$ $x_{k_{\mathbf{a},\alpha}} x_{l_{\mathbf{a},\alpha}} t_{q_{\mathbf{a},\alpha}}$ is a target monomial.  If there exists at least one such arrow $\mathbf{a}$ and one such triple $(t_{q_{\mathbf{a},1}},$ $k_{\mathbf{a},1},$ $l_{\mathbf{a},1})$, then $\Pi$ is spinal; indeed, the corollary implies that $\rho^{k_{\mathbf{a},1}l_{\mathbf{a},1}}_{pq_{\mathbf{a},1}}$ is in the spine of $T_{\Pi, k_{\mathbf{a},1}}$ and $T_{\Pi, l_{\mathbf{a},1}}$.  Having identified an index $k$ such that $T_{\Pi,k}$ has non-empty spine, one can compute the full spine by finding, for each arrow $\mathbf{a}$, the set of $\alpha$ for which $\{ k_{\mathbf{a},\alpha}, l_{\mathbf{a},\alpha} \}$ $=$ $\{ k, l'_{\mathbf{a},\alpha} \}$ with $d_{l'_{\mathbf{a},\alpha}} > 0$, and including $\rho^{kl'_{\mathbf{a},\alpha}}_{pq_{\mathbf{a},\alpha}}$ (resp.\ $\rho^{l'_{\mathbf{a},\alpha}k}_{pq_{\mathbf{a},\alpha}}$) in the spine if $k < l'_{\mathbf{a},\alpha}$ (resp.\ $l'_{\mathbf{a},\alpha} < k$). 
\end{rem}

\begin{exmp} \label{exmp:arrProc} Consider the case 
\[
  \begin{array}{rcl}
    \script{O} & = & \{1,\, x_1,\, x_2,\, x_3  \}\ \subseteq\ \mathbb{T}^3,\vspace{.05in}\\
    \partial \script{O} & = & \{x_1^2,\, x_1 x_2,\, x_1 x_3,\, x_2^2,\, x_2 x_3,\, x_3^2 \},\vspace{.05in}\\
    \Pi & = & \left< 1, 2, 3 \right>,\ \ \mdeg(\Pi) = (1,1,1)^{tr},
  \end{array}
\]
and index the elements of $\script{O}$ and $\partial \script{O}$ in the order displayed.  We will use the technique of Remark \ref{rem:SpAr} to show that $\Pi$ is spinal and to compute the spine of the trace syzygy $T_{\Pi, 1}$.  To begin, we note that $m$ $=$ $x_1 x_2 x_3$ is a target monomial in three ways:
\[
    x_1 x_2 x_3\ =\ x_2 x_3 t_2\ =\ x_1 x_3 t_3\ =\ x_1 x_2 t_4,
\]
moreover, the (ordinary) degree of any target monomial in this case can be at most $3$.  Since an arrow $t_p$ $\rightarrow$ $m'$ with displacement vector $(1,1,1)^{tr}$ will have a head of degree 4 if $p$ $=$ $2, 3, \text{ or } 4$, we see that $t_1$ $\rightarrow$ $m$ is the only possible arrow with this displacement vector that can be drawn from a monomial in $\script{O}$ to a target monomial.  Furthermore, $m$ is a $(1,2)$-target monomial, and $d_2$ $=$ $1 > 0$.  From this we conclude that $\Pi$ is spinal, since $T_{\Pi,1}$ has at least $\rho^{1,2}_{1,4}$ in its spine.  To compute the full spine of $T_{\Pi,1}$, we examine the witness-pairs $(k,l)$ for $m$ that include $1$ as a member, and note that $m$ is also a $(1,3)$-target monomial and $d_3$ $=$ $1 > 0$; whence, the spine of $T_{\Pi,1}$ has two members: $\rho^{1,2}_{1,4}$ and $\rho^{1,3}_{1,3}$.  By symmetry, similar results hold for $T_{\Pi,2}$ and $T_{\Pi,3}$.
\end{exmp}

We will conclude this paper with three extended examples of trace syzygies (in Sections \ref{sec:TSyzEx1}, \ref{sec:TSyzEx2}, and \ref{sec:2varGenEx}).  The first of these computes the trace syzygies with non-empty spines when $\script{O}$ $=$ $\{ 1, x_1, x_2 \}$ $\subseteq$ $\mathbb{T}^2$.  Using these syzygies, one demonstrates that the $\script{O}$-border basis scheme is isomorphic to six-dimensional affine space, a result that previously appeared as an example in \cite[Sec.\ 5.2]{Huib:HaimanUmu}.  The second example computes the trace syzygies with non-empty spines when $\script{O}$ $=$ $\{1, x_1\}$ $\subseteq$ $\mathbb{T}^3$.  Using these syzygies and the Jacobi identity syzygies previously computed in Example \ref{exmp:JacId3}, one can show that the $\script{O}$-border basis scheme is isomorphic to six-dimensional affine space; as noted in the introduction, this is a special case of \cite[Cor.\ 3.13]{Robbiano:BorderAndGrobner}.  The third example considers an arbitrary order ideal $\script{O}$ $\subseteq$ $\mathbb{T}^2$; using the trace syzygies with non-empty spines, and assuming that $\operatorname{char}(\gf)$ $=$ $0$, it is shown that the ideal of the $\script{O}$-border basis scheme $\mathbb{B}_{\script{O}}$ is generated by a subset of the $\rho^{1,2}_{pq}$ of size equal to the codimension of $\mathbb{B}_{\script{O}}$ in $\Spec(\gf[\mathbf{c}])$.



\section{Trace syzygies for $\script{O}$ $=$ $\{1, x_1, x_2\}$ $\subseteq$ $\mathbb{T}^2$} \label{sec:TSyzEx1}

Our goal is to compute the trace syzygies that have non-empty spines in this case.  Note first of all that the border is 
\[
    \partial \script{O}\ =\ \{ x_1^2,\, x_1 x_2,\, x_2^2 \};
\]
as usual, we will index the elements of $\script{O}$ and $\partial \script{O}$ as displayed.

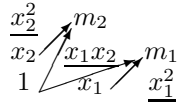
\begin{figure}
\begin{picture}(300,60)
    \put(150,0){\makebox(0,0)[bc]{$
       \begin{array}{ccc}
          \underline{x_2^2} & m_2 & {}    \\
           x_2  & \underline{x_1 x_2} & m_1  \\
           1 & x_1 & \underline{x_1^2} 
    \end{array}$}}
    \put(156,5){\vector(1,1){12}}
    \put(129,18){\vector(1,1){12}}
    \put(129,4){\vector(1,2){12}}
    \put(129,4){\vector(3,1){36}}
\end{picture}
 \caption{Diagram of $\script{O}$ $=$ $\{1, x_1, x_2\}$ showing the four arrows with tail in $\script{O}$, head a target monomial, and a $(1,2)$-good direction vector.  The border monomials are underlined and the target monomials are $m_1$ $=$ $x_1^2 x_2$ and $m_2$ $=$ $x_1 x_2^2$.} \label{fig:1x1x2}
\end{figure}

To find the spinal multi-degrees $(d_1, d_2)^{tr}$, we use the idea of Remark \ref{rem:SpAr}, so we first determine the $(1,2)$-target monomials; by inspection, there are two:
\[
    m_1 = x_1^2 x_2,\ \text{and } m_2 = x_1 x_2^2.
\]  
Next, we find all the arrows $\mathbf{a}$ $=$ $(t_p \rightarrow m)$ such that $t_p$ $\in$ $\script{O}$, $m$ $\in$ $\{ m_1, m_2 \}$, and the direction vector
\[
  (d_1,d_2)^{tr}\ =\ \mdeg(m) - \mdeg(t_p)
\]
is $(1,2)$-good, which in this case means that both components are positive.  By inspection, there are four such arrows, two having direction vector $(1,1)^{tr}$, one having direction vector $(1,2)^{tr}$, and one having direction vector $(2,1)^{tr}$; see Figure \ref{fig:1x1x2}.  Each of these arrows demonstrates that its corresponding direction vector is spinal, as the criteria given in Remark \ref{rem:SpAr} are easily verified. We proceed to compute the trace syzygies associated to each of these multi-degrees.

The multiplication matrices are
\[
   \gfmm_1 = \left(
\begin{array}{ccc}
 0 & c_{1,1} & c_{1,2} \\
 1 & c_{2,1} & c_{2,2} \\
 0 & c_{3,1} & c_{3,2}
\end{array}
\right),\ \ \gfmm_2 = \left(
\begin{array}{ccc}
 0 & c_{1,2} & c_{1,3} \\
 0 & c_{2,2} & c_{2,3} \\
 1 & c_{3,2} & c_{3,3}
\end{array}
\right),
\]
and the commutator matrix is 
\[
   \begin{array}{l}
   \hspace*{1in} [\gfmm_1, \gfmm_2]\ = \ \left(
\begin{array}{ccc}
 \rho ^{1,2}_{1,1} & \rho ^{1,2}_{1,2} & \rho ^{1,2}_{1,3}\vspace{.05in} \\
 \rho ^{1,2}_{2,1} & \rho ^{1,2}_{2,2} & \rho ^{1,2}_{2,3}\vspace{.05in} \\
 \rho ^{1,2}_{3,1} & \rho ^{1,2}_{3,2} & \rho ^{1,2}_{3,3}
\end{array}
\right)\ = \vspace{.1in}\\
\small{
\left(
\begin{array}{ccc}
 0 & -c_{1,2} c_{2,1}+c_{1,1} c_{2,2}-c_{1,3}
   c_{3,1}+c_{1,2} c_{3,2} & -c_{1,2}
   c_{2,2}+c_{1,1} c_{2,3}-c_{1,3}
   c_{3,2}+c_{1,2} c_{3,3} \\
 0 & c_{1,2}-c_{2,3} c_{3,1}+c_{2,2} c_{3,2} &
   -c_{2,2}^2+c_{3,3} c_{2,2}+c_{1,3}+c_{2,1}
   c_{2,3}-c_{2,3} c_{3,2} \\
 0 & c_{3,2}^2-c_{2,1} c_{3,2}-c_{1,1}+c_{2,2}
   c_{3,1}-c_{3,1} c_{3,3} & -c_{1,2}+c_{2,3}
   c_{3,1}-c_{2,2} c_{3,2}
\end{array}
\right).
}
   \end{array} 
\]

\begin{exmp} \label{exmp:TSyz1} $\mathbf{ \Pi = \left<1, 2\right>,\ \mdeg(\Pi) = (1,1)^{tr}.}$   The expression defining the trace syzygy $T_{\Pi,1}$ is 
\[ 
    \Tr \left( [\gfmm_1, \gfmm_2]  \right)\ =\ \rho^{1,2}_{1,1} + \rho^{1,2}_{2,2} + \rho^{1,2}_{3,3};
\]
noting that $\rho^{1,2}_{1,1}$ $=$ $0$ (it is ``trivially $0$''), we see that 
\begin{equation} \label{eqn:Treln1}
  T_{\Pi, 1}\ =\  \left( \rho^{1,2}_{2,2} + \rho^{1,2}_{3,3}\,  =\, 0 \right)\ \Rightarrow\ \rho^{1,2}_{3,3} = - \rho^{1,2}_{2,2}. 
\end{equation}
The two $\rho$'s in the spine correspond to the arrows with displacement vector $(1,1)^{tr}$ in Figure \ref{fig:1x1x2}.  One sees easily that $T_{\Pi,2}$ $=$ $-T_{\Pi,1}$ in this case.
\end{exmp}

\begin{exmp} \label{exmp:TSyz2} $\mathbf{\Pi = \left< 1,1,2 \right>, \mdeg(\Pi) = (2,1)^{tr}}$.  The expression defining the syzygy $T_{\Pi,1}$ is
\[
  \begin{array}{rcl}
    \Tr\left( [\gfmm_1,\gfmm_1] \cdot \gfmm_2 + \gfmm_1 \cdot [\gfmm_1, \gfmm_2] \right) & = & \Tr \left( \gfmm_1 \cdot [\gfmm_1, \gfmm_2] \right) \vspace{.1in}\\
    {} & = &  
  \left\{\begin{array}{l}
            c_{1,1}\, \rho^{1,2}_{2,1} + c_{2,1}\, \rho
   ^{1,2}_{2,2} + c_{3,1}\, \rho^{1,2}_{2,3} \vspace{.05in}\ + \\ c_{1,2}\, \rho
   ^{1,2}_{3,1} + c_{2,2}\, \rho^{1,2}_{3,2} + c_{3,2}\, \rho
   ^{1,2}_{3,3} + \rho^{1,2}_{1,2}
         \end{array} \right\};
  \end{array}
\]
removing the terms involving $\rho^{1,2}_{2,1}$ and $\rho^{1,2}_{3,1}$, which are trivially $0$, we obtain
\[
  T_{\Pi, 1}\ =\ \left( c_{2,1}\, \rho
   ^{1,2}_{2,2} + c_{3,1}\, \rho^{1,2}_{2,3} + c_{2,2}\, \rho^{1,2}_{3,2} + c_{3,2}\, \rho
   ^{1,2}_{3,3} + \rho^{1,2}_{1,2}  = 0 \right).
\]
We can further simplify by replacing $\rho^{1,2}_{3,3}$ by $-\rho^{1,2}_{2,2}$ and regrouping; the result is the relation
\begin{equation} \label{eqn:Treln2}
    (c_{2,1} - c_{3,2})\, \rho
   ^{1,2}_{2,2} + c_{3,1}\, \rho^{1,2}_{2,3} + c_{2,2}\, \rho^{1,2}_{3,2}  + \rho^{1,2}_{1,2}\ =\ 0.
\end{equation}
In this case there is only one element in the spine, which corresponds to the single arrow with direction vector $(2,1)^{tr}$ in Figure \ref{fig:1x1x2}.

 The syzygy $T_{\Pi,2}$ is defined by the expression
\[
    \Tr \left( [\gfmm_2,\gfmm_1]\cdot \gfmm_1 + \gfmm_1 \cdot [\gfmm_2, \gfmm_1] \right)\ = \ -2 \Tr \left( \gfmm_1 \cdot [\gfmm_1,\gfmm_2]\right);
\]
whence, $T_{\Pi,2}$ $=$ $-2 \cdot T_{\Pi,1}$ in this case, which in particular exemplifies Corollary \ref{cor:T3}.
\end{exmp}

\begin{exmp} \label{exmp:TSyz3} $\mathbf{\Pi = \left<1, 2, 2\right>,\ \mdeg(\Pi) = (1,2)^{tr}}$.  The expression defining $T_{\Pi,2}$ is 
\[
  \begin{array}{rcl}
  \Tr([\gfmm_2,\gfmm_1]\cdot \gfmm_2 + \gfmm_1 \cdot [\gfmm_2, \gfmm_2]) & =  & \Tr([\gfmm_2, \gfmm_1] \cdot \gfmm_2)\vspace{.05in}\\ 
 {} & = & -\Tr([\gfmm_1, \gfmm_2] \cdot \gfmm_2)\vspace{.05in}\\
 {} & = &  
  - \left\{   \begin{array}{l}
 c_{1,2}\, \rho^{1,2}_{2,1} +  c_{2,2}\, \rho^{1,2}_{2,2} + 
   c_{3,2}\, \rho^{1,2}_{2,3}\ +\vspace{.05in}\\  c_{1,3}\, \rho
   ^{1,2}_{3,1} +
  c_{2,3}\, \rho^{1,2}_{3,2} +  c_{3,3}\,
   \rho^{1,2}_{3,3} +  \rho^{1,2}_{1,3}
             \end{array} \right\} ;
  \end{array}
\]
proceeding as in the $(2,1)$-case, one obtains the relation
\begin{equation} \label{eqn:Treln3}
      (c_{2,2} - c_{3,3}) \, \rho^{1,2}_{2,2} + 
   c_{3,2}\, \rho^{1,2}_{2,3} +  
  c_{2,3}\, \rho^{1,2}_{3,2} +  \rho^{1,2}_{1,3}\ =\ 0
\end{equation}
and checks that $T_{\Pi,1}$ $=$ $-2 \cdot T_{\Pi,2}$.
\end{exmp}

\begin{rem} \label{rem:UmuRef}

The relations (\ref{eqn:Treln1}), (\ref{eqn:Treln2}), and (\ref{eqn:Treln3}) were presented in  \cite[Sec.\ 5.2]{Huib:HaimanUmu}, but in a different notation, and with a completely different derivation (given in \cite[Sec.\ 6.2]{Huib:HaimanUmu}) that is far less general than the approach used here.

\end{rem}

\begin{rem} \label{rem:compIntRem}

The $\script{O}$-border basis scheme $\mathbb{B}_{\script{O}}$ studied in the preceding examples is a subvariety of $9$-dimensional affine space having coordinate functions 
\[
    \{ c_{ij}\, \mid\, 1 \leq i,j \leq  3 \},
\]
and cut out by the six non-zero polynomials
\[
     \{ \rho^{1,2}_{ij}\, \mid\, 1 \leq i \leq 3,\ 2 \leq j \leq 3 \}.
\]
The relations (\ref{eqn:Treln1}), (\ref{eqn:Treln2}), and (\ref{eqn:Treln3}) show that 
\begin{equation} \label{eqn:threeRhos}
   \rho^{1,2}_{3,3},\ \rho^{1,2}_{1,2},\ \rho^{1,2}_{1,3}\ \in 
\left( \rho^{1,2}_{2,2},\ \rho^{1,2}_{2,3},\ \rho^{1,2}_{3,2} \right);
\end{equation}
so in fact the ideal of $\mathbb{B}_{\script{O}}$ is generated by three polynomials.  Furthermore, as noted in \cite[Sec.\ 5.2]{Huib:HaimanUmu}, the three ideal generators express $c_{1,1}$, $c_{1,2}$, and $c_{1,3}$ in terms of the six remaining $c$'s, which implies that 
\[
  \mathbb{B}_{\script{O}}\ =\ \Spec(\gf[c_{2,1},\, c_{2,2},\, c_{2,3},\, c_{3,1},\, c_{3,2},\, c_{3,3}])
\]
is an affine space of dimension six.  In particular, the trace syzygies allow us in this case to discard certain $\rho$'s and demonstrate that $\mathbb{B}_{\script{O}}$ is a complete intersection in $\Spec(\gf[\mathbf{c}])$.  This will be generalized to any two-variable border basis scheme in Section \ref{sec:2varGenEx} (in characteristic $0$). 
\end{rem}

\section{Trace syzygies for $\script{O}$ $=$ $\{1, x_1 \}$ $\subseteq$ $\mathbb{T}^3$} \label{sec:TSyzEx2}

Recall that we computed the Jacobi identity syzygies associated to this order ideal in Example \ref{exmp:JacId3}, where the multiplication matrices (\ref{eqn:JExMultMats}) and their commutators (\ref{eqn:JExComms}) can be found. (Note that all of the $\rho^{kl}_{pq}$ are non-zero in this case.)  By multiplying the monomials $t_1$ $=$ $1$ and $t_2$ $=$ $x_1$ by $x_1 x_2$, $x_1 x_3$, and $x_2 x_3$, we find that there are six target monomials:
\[
  \begin{array}{l}
      (1,2)\text{-target monomials}:\ x_1 x_2,\ x_1^2 x_2\\
      (1,3)\text{-target monomials}:\ x_1 x_3,\ x_1^2 x_3\\
      (2,3)\text{-target monomials}:\ x_2 x_3,\ x_1 x_2 x_3.
  \end{array}
\]
To find the spinal multi-degrees $(d_1, d_2, d_3)^{tr}$, we again use the method of Remark \ref{rem:SpAr}: we begin by computing the vectors $\mdeg(m)$ $-$ $\mdeg(t_p)$ as $m$ runs through the target monomials and $t_p$ runs through the monomials in $\script{O}$, and retaining those that are good.  This results in the following list:
\[
    (d_1, d_2, d_3)^{tr}\ =\ (1,1,0)^{tr},\, (1,0,1)^{tr},\, (0,1,1)^{tr},\, (2,1,0)^{tr},\, (2,0,1)^{tr},\, (1,1,1)^{tr}.
\]

There are two arrows $t_p$ $\rightarrow$ $m$ associated to each of the first three of these multi-degrees, and only one to each of the last three.  One checks easily that each of these arrows demonstrates (as in Remark \ref{rem:SpAr}) that its associated multi-degree is spinal.

We now exhibit a trace syzygy with non-empty spine associated to each of the six spinal multi-degrees.  In each case, we display the expression (\ref{eqn:TSyzSum}) that defines the trace syzygy, followed by the trace syzygy itself.
\begin{description}
  \item[$\mathbf{\Pi = \left<1, 2\right>,\ \mdeg(\Pi) = (1,1,0)^{tr}}$]  
\[
   \text{Expression: } \Tr \left( [\gfmm_1, \gfmm_2] \right),\ \ T_{\Pi,1}\ =\ \left( \rho^{1,2}_{1,1} + \rho^{1,2}_{2,2}\, =\, 0 \right)
\] 
  \item[$\mathbf{\Pi = \left<1, 3\right>,\ \mdeg(\Pi) = (1,0,1)^{tr}}$]
\[
    \text{Expression: } \Tr \left( [\gfmm_1, \gfmm_3] \right),\ \ T_{\Pi,1}\ =\  \left( \rho^{1,3}_{1,1} + \rho^{1,3}_{2,2}\, =\, 0 \right)
\] 
  \item[$\mathbf{\Pi = \left<2, 3\right>,\ \mdeg(\Pi) = (0,1,1)^{tr} }$]
\[
    \text{Expression: } \Tr \left( [\gfmm_2, \gfmm_3] \right),\ \ T_{\Pi,2} \ =\ \left( \rho^{2,3}_{1,1} + \rho^{2,3}_{2,2}\, =\, 0 \right) 
\]
  \item[$\mathbf{\Pi = \left<1, 1, 2\right>,\ \mdeg(\Pi) = (2,1,0)^{tr}}$]
\[
  \begin{array}{c} 
      \text{Expression: } \Tr \left( [\gfmm_1, \gfmm_1]\cdot \gfmm_2 + \gfmm_1\cdot[\gfmm_1, \gfmm_2] \right), \vspace{.05in}\\
      T_{\Pi, 1}\ =\  \left( c_{1,3}\, \rho^{1,2}_{2,1} + c_{2,3}\, \rho^{1,2}_{2,2} + \rho
  ^{1,2}_{1,2}\, =\, 0 \right)
  \end{array}
\]
  \item[$\mathbf{\Pi = \left<1, 1, 3\right>,\ \mdeg(\Pi) = (2,0,1)^{tr}}$]
\[
  \begin{array}{c} 
      \text{Expression: } \Tr \left( [\gfmm_1, \gfmm_1]\cdot \gfmm_3 + \gfmm_1\cdot[\gfmm_1, \gfmm_3] \right), \vspace{.05in}\\
     T_{\Pi,1}\ =\ \left(c_{1,3}\, \rho^{1,3}_{2,1} + c_{2,3}\, \rho^{1,3}_{2,2} + \rho
  ^{1,3}_{1,2}\, =\, 0 \right)
  \end{array}
\] 
  \item[$\mathbf{\Pi = \left<1, 2, 3\right>,\ \mdeg(\Pi) = (1,1,1)^{tr}}$]  Note that since $x_1 x_2 x_3$ is a $(2,3)$-target monomial (only), we must use either $2$ or $3$ as distinguished index to obtain a non-empty spine; we choose $2$.
\[
    \begin{array}{c} 
       
 \text{Expression: } \Tr \left( -[\gfmm_1, \gfmm_2]\cdot \gfmm_3 + \gfmm_1\cdot[\gfmm_2, \gfmm_3] \right), \vspace{.05in}\\
      T_{\Pi, 2}\ =\  \left(
          \left\{ \begin{array}{l}-c_{1,2}\, \rho^{1,2}_{1,1} - c_{2,2}\, \rho^{1,2}_{1,2} - c_{1,5}\, \rho
  ^{1,2}_{2,1} - c_{2,5}\, \rho^{1,2}_{2,2}\, + \vspace{.05in}\\
     c_{1,3}\, \rho
  ^{2,3}_{2,1} + c_{2,3}\, \rho^{2,3}_{2,2} + \rho^{2,3}_{1,2} 
          \end{array}\right\}  = 0 \right)
  \end{array}
\] 
\end{description}   
\medskip

We now use these syzygies to confirm that the $\script{O}$-border basis scheme $\mathbb{B}_{\script{O}}$ is isomorphic to an affine space. This is a special case of \cite[Cor.\ 3.13]{Robbiano:BorderAndGrobner}, which asserts that for any order ideal of the form 
\[
    \script{O}\ =\ \{ 1, x_k, x_k^2, \dots, x_k^{\mu-1} \} \subseteq \mathbb{T}^n,
\]
one has that $\mathbb{B}_{\script{O}}$ is isomorphic to the affine space $\mathbb{A}^{\mu \cdot n}$; consequently, for $\script{O}$ $=$ $\{1, x_1  \}$ $\subseteq$ $\mathbb{T}^3$, we have that $\mathbf{B}_{\script{O}}$ is isomorphic to an affine space of dimension $2 \cdot 3$ $=$ $6$.  In fact, we will show that
\begin{equation} \label{eqn:1x1CoordFuncs}
    \mathbb{B}_{\script{O}}\ =\ \Spec(\gf[\mathbf{c}]/I(\mathbb{B}_{\script{O}}))\ =\ \Spec(\gf[c_{1,1},\, c_{1,2},\, c_{1,3},\, c_{2,1},\, c_{2,2},\, c_{2,3}]).
\end{equation}

To do this, we first note that the following generators of $I(\mathbb{B}_{\script{O}})$:
\begin{equation} \label{eqn:fourRhos}
    \begin{array}{rcl} 
	\rho^{1,2}_{1,1} & = & c_{1,3}\, c_{2,1}-c_{1,4},\vspace{.05in}\\
	\rho^{1,2}_{2,1} & = & c_{1,1}+c_{2,1}\, c_{2,3}-c_{2,4},\vspace{.05in}\\
	\rho^{1,3}_{1,1} & = &  c_{1,3}\, c_{2,2}-c_{1,5},\vspace{.05in}\\
	\rho^{1,3}_{2,1} & = &  c_{1,2}+c_{2,2}\, c_{2,3}-c_{2,5},
    \end{array}
\end{equation}
demonstrate that each of $c_{1,4}$, $c_{2,4}$, $c_{1,5}$, $c_{2,5}$, respectively, is congruent modulo $I(\mathbb{B}_{\script{O}})$ to a polynomial in the six $c$'s that generate the polynomial ring in (\ref{eqn:1x1CoordFuncs}).  We now have
\begin{lem} \label{lem:1x1Claim} The eight remaining generators of $I(\mathbb{B}_{\script{O}})$, namely
\[
    \rho^{1,2}_{1,2},\ \rho^{1,2}_{2,2},\ \rho^{1,3}_{1,2},\ \rho^{1,3}_{2,2},\ \rho^{2,3}_{1,1},\ \rho^{2,3}_{1,2},\ \rho^{2,3}_{2,1},\ \rho^{2,3}_{2,2},
\]
all lie in the ideal 
\[
  J\ =\ \left( \rho^{1,2}_{1,1},\ \rho^{1,2}_{2,1},\ \rho^{1,3}_{1,1},\ \rho^{1,3}_{2,1} \right)  \subseteq \Zc
\]
 generated by the four $\rho$'s in \emph{(\ref{eqn:fourRhos})}.
\end{lem}

Assuming the Lemma, we see that 
\[
   \gf[\mathbf{c}]/I(\mathbb{B}_{\script{O}}))\ =\ \gf[\mathbf{c}]/(\rho^{1,2}_{1,1},\, \rho^{1,2}_{2,1},\, \rho^{1,3}_{1,1},\, \rho^{1,3}_{2,1})\ =\ \gf[c_{1,1},\, c_{1,2},\, c_{1,3},\, c_{2,1},\, c_{2,2},\, c_{2,3}],
\]
as desired, so it suffices to prove the Lemma.
\medskip

\emph{Proof of Lemma \ref{lem:1x1Claim}}: Several of the cases are easy to check:
\[
    \begin{array}{l}
        \rho^{1,2}_{2,2}\ =\ -\rho^{1,2}_{1,1} \in  J\ \text{by } T_{\left<1,2\right>,1},\vspace{.05in}\\
	\rho^{1,3}_{2,2}\ =\ -\rho^{1,3}_{1,1} \in J\  \text{by } T_{\left<1,3\right>,1},\vspace{.05in}\\
	\rho^{1,2}_{1,2}\ =\ -c_{1,3}\, \rho^{1,2}_{2,1} - c_{2,3}\, \rho^{1,2}_{2,2}\ =\ -c_{1,3}\, \rho^{1,2}_{2,1} + c_{2,3}\, \rho^{1,2}_{1,1} \in J,\ \text{by } T_{\left<1,1,2\right>,1},\vspace{.05in}\\
	\rho^{1,3}_{1,2}\ =\ -c_{1,3}\, \rho^{1,3}_{2,1} - c_{2,3}\, \rho^{1,3}_{2,2}\ =\ -c_{1,3}\, \rho^{1,3}_{2,1} + c_{2,3}\, \rho^{1,3}_{1,1} \in J\ \text{ by } T_{\left<1,1,3\right>,1}.
    \end{array}
\]
Furthermore, $T_{\left<2,3\right>,2}$ implies that $\rho^{2,3}_{2,2}$ $=$ $-\rho^{2,3}_{1,1}$, so it remains to establish the following conclusions:
\[
    \rho^{2,3}_{1,1} \in J,\ \rho^{2,3}_{1,2} \in J,\ \text{and } \rho^{2,3}_{2,1} \in J.
\]
In light of the foregoing, one has that $T_{\left<1,2,3\right>,2}$ implies that
\[
    \left( \rho^{2,3}_{1,1} \in J \text{ and } \rho^{2,3}_{2,1} \in J \right) \Rightarrow \rho^{2,3}_{1,2} \in J,
\]
but the trace syzygies appear inadequate to establish the two inclusions in the antecedent.  To do this, we use the Jacobi identity syzygy (\ref{eqn:JsyzEg3list}):
\[
    J^{1,2,3}_{2,1}\  =\  \left( \left\{ \begin{array}{l} c_{2,2}\, \rho^{1,2}_{1,1} +(c_{2,5} - c_{1,2})\, \rho
   ^{1,2}_{2,1} - c_{2,2}\, \rho
   ^{1,2}_{2,2} - c_{2,1}\, \rho^{1,3}_{1,1}\ + \vspace{.05in}\\
   (c_{1,1} - c_{2,4})\, \rho
   ^{1,3}_{2,1} + c_{2,1}\, \rho
   ^{1,3}_{2,2} + c_{2,3}\, \rho^{2,3}_{2,1} + \rho
   ^{2,3}_{1,1} - \rho^{2,3}_{2,2} \end{array}\right\}  = 0 \right).
\]
By replacing $\rho^{k,l}_{2,2}$ with $-\rho^{k,l}_{1,1}$ and regrouping, we obtain the relation
\[
    -2\rho^{2,3}_{1,1} -2 c_{2,2}\, \rho^{1,2}_{1,1} +2 c_{2,1}\, \rho^{1,3}_{1,1}\ =\  (c_{2,5} - c_{1,2})\, \rho
   ^{1,2}_{2,1}  + (c_{1,1} - c_{2,4})\, \rho
   ^{1,3}_{2,1}\ + c_{2,3}\, \rho^{2,3}_{2,1}.
\]
By inspection, none of the $c$'s on the RHS appear in any of the $\rho$'s on the LHS (see (\ref{eqn:JExComms})); this implies that both sides of the equation must be $0$ as elements of $\Zc$, which yields 
\[
    \begin{array}{rcl}
\rho^{2,3}_{1,1} & = &  - c_{2,2}\, \rho^{1,2}_{1,1} +  c_{2,1}\, \rho^{1,3}_{1,1} \in J,\ \text{and}\vspace{.05in}\\
c_{2,3}\, \rho^{2,3}_{2,1} & = & (c_{1,2} - c_{2,5})\, \rho
   ^{1,2}_{2,1}  - (c_{1,1} - c_{2,4})\, \rho
   ^{1,3}_{2,1}\vspace{.05in}\\
{} & = & (\rho^{1,3}_{2,1} - c_{2,2}\,c_{2,3})\, \rho
   ^{1,2}_{2,1}  - (\rho^{1,2}_{2,1} - c_{2,1}\,c_{2,3})\, \rho
   ^{1,3}_{2,1}\vspace{.05in}\\
{} & = & c_{2,1}\,c_{2,3}\, \rho
   ^{1,3}_{2,1} - c_{2,2}\,c_{2,3}\, \rho
   ^{1,2}_{2,1}\vspace{.05in}\\
{} & \Rightarrow & \rho^{2,3}_{2,1}\ =\ c_{2,1}\, \rho
   ^{1,3}_{2,1} - c_{2,2}\, \rho
   ^{1,2}_{2,1} \in J.
    \end{array}
\]
This completes the proof of the lemma and of the assertion (\ref{eqn:1x1CoordFuncs}). \qedsymbol

\section{Trace syzygies for an arbitrary $\script{O}$ $\subseteq$ $\mathbb{T}^2$} \label{sec:2varGenEx}

In this final section of the paper, we assume that $\operatorname{char}(\gf)$ $=$ $0$.  Let $\script{O}$ $\subseteq$ $\mathbb{T}^2$ be an arbitrary order ideal of size $\mu$ and with border $\partial \script{O}$ of size $\nu$.  We will use trace syzygies to show that the border basis scheme 
\[
    \mathbb{B}_{\script{O}} = \Spec(\gf[\mathbf{c}]/I(\mathbb{B}_{\script{O}}))\ \subseteq\ \mathbb{A}^{\mu \nu} ,\ \ I(\mathbb{B}_{\script{O}}) = (\rho^{1,2}_{pq}),\ 1 \leq p, q \leq \mu,
\]
is cut out by $(\nu - 2)\mu$ of the (not-trivially-$0$) $\rho^{1,2}_{pq}$.  

To clarify this result, recall that $\mathbb{B}_{\script{O}}$ is irreducible and nonsingular of dimension $2\mu$, since it is an open subscheme of the Hilbert scheme $\Hilb^{\mu}_{\mathbb{A}^2}$, which is well-known to have these properties. (This was first proved by Fogarty \cite[Th.\ 2.4]{Fogarty:AFAS}; see \cite[Prop.\ 2.4]{Haiman:CN-HS} for a proof making use of the concept, if not the terminology, of border basis schemes.)  Since the number of indeterminates $c_{ij}$ is equal to $\mu \nu$, our result says that $\mathbb{B}_{\script{O}}$ is a complete intersection in $\mathbb{A}^{\mu \nu}$, and identifies an explicit minimal generating set of the ideal $I(\mathbb{B}_{\script{O}})$.

\begin{figure}
\begin{picture}(300,132)
    \put(150,0){\makebox(0,0)[bc]{$
       \begin{array}{ccccccccccc}
          \bullet & \bullet & \bullet &{}&{}&{}&{}&{}&{}&{}&{}\\
          * & * & * & \bullet &{}&{}&{}&{}&{}&{}&{}\\
          . & . & * & \bullet & \bullet & \bullet & \bullet &{}&{}&{}&{}\\
          . & . & . & * & * & * & * & \bullet & {}&{}&{}\\
          . & . & . & . & . & . & * & \bullet & {}&{}&{}\\
          . & . & . & . & . & . & * & \bullet & {}&{}&{}\\
          . & . & . & . & . & . & * & \bullet & {}&{}&{}\\
          . & . & . & . & . & . & * & \bullet & \bullet & \bullet &{}\\
          . & . & . & . & . & . & . & * & * & * & \bullet\\
          x_2 & . & . & . & . & . & . & . & . & * & \bullet\\
          1 & x_1 & . & . & . & . & . & . & . & * & \bullet\\
       \end{array}$}}
    \put(66,0){\line(0,1){120}}
    \put(66,120){\line(1,0){48}}
    \put(114,120){\line(0,-1){24}}
    \put(114,96){\line(1,0){60}}
    \put(174,96){\line(0,-1){60}}
    \put(174,36){\line(1,0){48}}
    \put(222,36){\line(0,-1){36}}
    \put(222,0){\line(-1,0){156}}
    \put(72,117){\line(0,1){6}}
    \put(93,117){\line(0,1){6}}
    \put(110,117){\line(0,1){6}}
    \put(113,103){\line(1,1){10}}
    \put(125,93){\line(0,1){6}}
    \put(140,93){\line(0,1){6}}
    \put(155,93){\line(0,1){6}}
    \put(170,93){\line(0,1){6}}
    \put(172,78){\line(1,1){10}}
    \put(172,66){\line(1,1){10}}
    \put(172,54){\line(1,1){10}}
    \put(172,42){\line(1,1){10}}
    \put(185,33){\line(0,1){6}}
    \put(200,33){\line(0,1){6}}
    \put(215,33){\line(0,1){6}}
    \put(217,18){\line(1,1){10}}
    \put(217,6){\line(1,1){10}}
\end{picture}
 \caption{Diagram of an order ideal $\script{O}$ $\subseteq$ $\mathbb{T}^2$, with the monomials in $\script{O}$ enclosed in a box.  The exposable monomials are shown as $*$'s, and the border monomials (those in $\partial \script{O}$) are shown as $\bullet$'s.  The correspondence indicated by the line segments shows that the number of $*$'s is one less than the number of $\bullet$'s (of which there are $\nu$).} \label{fig:2varcase}
\end{figure}
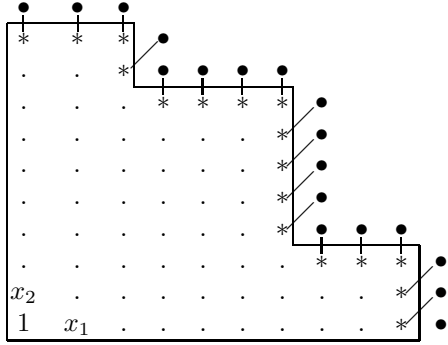

Consider first Figure \ref{fig:2varcase}, which shows the diagram of a typical order ideal in two variables, along with its border. We will call a monomial $t_q$ $\in$ $\script{O}$ \textbf{exposable} if either $x_1 t_q$ $\notin$ $\script{O}$ or $x_2 t_q$ $\notin$ $\script{O}$ (in the figure, the exposable monomials are shown as $*$'s). By the analysis in Section \ref{sec:BBSidGens}, we know that $\rho^{1,2}_{pq}$ is not trivially $0$ if and only if $t_q$ is exposable (since then one is in either Case 3 or Case 4).  Therefore: 

\begin{lem} \label{lem:numRhos}
    The total number of $\rho^{1,2}_{pq}$ that are not trivially $0$ is equal to $(\nu - 1)\mu$. 
\end{lem}

\emph{Proof:}
    As indicated in Figure \ref{fig:2varcase}, the number of exposable monomials $t_q$ is equal to $\mu$ $-$ $1$, and for each such $q$, the value of $p$ in $\rho^{1,2}_{pq}$ can take on the values $1$ $\dots$ $\mu$.
 \qedsymbol
\medskip

Since there are only two variables $x_1$ and $x_2$, there are no nontrivial Jacobi identity syzygies among the $\rho$'s (Remark \ref{rem:noJacIn2Vars}).  On the other hand, there are plenty of trace syzygies:

\begin{lem} \label{lem:NumTSyz}
    The number of spinal multi-degrees $(d_1, d_2)^{tr}$ is equal to $\mu$, the number of monomials in $\script{O}$.
\end{lem}

\emph{Proof:}
     First of all, there are only two possible distinguished indices ($1$ and $2$), and the spines of $T_{\Pi, 1}$ and $T_{\Pi, 2}$  are equal, by Corollary \ref{cor:T3}.  Hence, to show that the $(1,2)$-good tuple $(d_1, d_2)^{tr}$ is spinal, it suffices to show that $T$ $=$ $T_{\Pi, 1}$ has non-empty spine, where $\Pi$ an ordered product such that $\mdeg(\Pi)$ $=$ $(d_1, d_2)^{tr}$.  By Corollary \ref{cor:TarMD}, $T$ has non-empty spine if and only if there is an arrow $\mathbf{a}$ $=$ $t_p$ $\rightarrow$ $m$ with direction vector $(d_1, d_2)^{tr}$, $t_p$ $\in$ $\script{O}$, and $m$ $=$ $x_1 x_2 t_q$ a $(1,2)$-target monomial; indeed, the existence of $\mathbf{a}$ $=$ $\arr(\rho^{1,2}_{pq})$ witnesses that $\rho^{1,2}_{pq}$ is in the spine of $T$.  

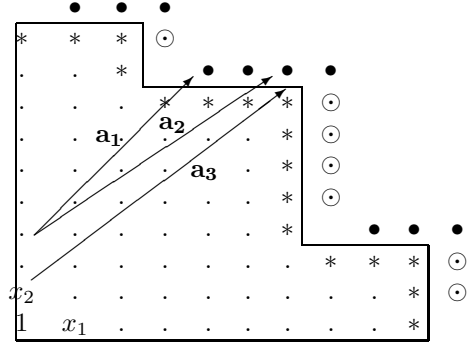
\begin{figure}
\begin{picture}(300,132)
    \put(150,0){\makebox(0,0)[bc]{$
       \begin{array}{ccccccccccc}
          {} & \bullet & \bullet & \bullet &{}&{}&{}&{}&{}&{}&{}\\
          * & * & * & \odot &{}&{}&{}&{}&{}&{}&{}\\
          . & . & * & {} & \bullet & \bullet & \bullet & \bullet&{}&{}&{}\\
          . & . & . & * & * & * & * & \odot & {}&{}&{}\\
          . & . & . & . & . & . & * & \odot & {}&{}&{}\\
          . & . & . & . & . & . & * & \odot & {}&{}&{}\\
          . & . & . & . & . & . & * & \odot & {}&{}&{}\\
          . & . & . & . & . & . & * & {} & \bullet & \bullet & \bullet\\
          . & . & . & . & . & . & . & * & * & * & \odot\\
          x_2 & . & . & . & . & . & . & . & . & * & \odot\\
          1 & x_1 & . & . & . & . & . & . & . & * & {}\\
       \end{array}$}}
    \put(66,0){\line(0,1){120}}
    \put(66,120){\line(1,0){48}}
    \put(114,120){\line(0,-1){24}}
    \put(114,96){\line(1,0){60}}
    \put(174,96){\line(0,-1){60}}
    \put(174,36){\line(1,0){48}}
    \put(222,36){\line(0,-1){36}}
    \put(222,0){\line(-1,0){156}}
    \put(73,40){\vector(1,1){60}}
    \put(73,40){\vector(3,2){90}}
    \put(72,23){\vector(4,3){96}}
    \put(96,75){$\mathbf{a_1}$}
    \put(120,80){$\mathbf{a_2}$}
    \put(132,62){$\mathbf{a_3}$}
\end{picture}
 \caption{In this diagram the $\bullet$'s represent the $(1,2)$-target monomials that can be targets of arrows in extreme position (\ie, are on ``horizontal steps''), and the $\odot$'s represent the remaining $(1,2)$-target monomials.  The $*$'s represent the exposable monomials $t_q$ $\in$ $\script{O}$; for each of these, $x_1 x_2 t_q$ is a $(1,2)$-target monomial.  Three arrows in extreme position are shown, having direction vectors $\mathbf{d_1}$ $=$ $(4,5)^{tr}$, $\mathbf{d_2}$ $=$ $(6,5)^{tr}$, and $\mathbf{d_3}$ $=$ $(6,7)^{tr}$, respectively.  In the ordering of such arrows defined in the text, we have $\mathbf{a_1}$ $\prec$ $\mathbf{a_2}$ $\prec$ $\mathbf{a_3}$. 
} \label{fig:2vExtVec}
\end{figure}

Given such an arrow $\mathbf{a}$, one can translate it to the ``northwest,'' keeping its source in $\script{O}$ and its target among the $(1,2)$-target monomials, until its source is along the left edge of $\script{O}$ (that is, the source is of the form $t_p$ $=$ $x_2^u$) and its target $m$ is on one of the ``horizontal steps'' in the diagram (more precisely, $m$ $=$ $x_1 x_2 t_q$ with $x_2 t_q$ $\notin$ $\script{O}$).  An arrow that has been translated as far as possible to the northwest as described will be said to be in \textbf{extreme} position, and the associated polynomial $\rho^{1,2}_{p,q}$, which lies in the spine of $T$, will be called an \textbf{extreme} $\rho$ --- see Figure \ref{fig:2vExtVec}.  

     At this point we know that $T$ has a non-empty spine if and only if there is an arrow $\mathbf{a}$ in extreme position and having $(1,2)$-good direction vector $\mdeg(\Pi)$ $=$ $(d_1, d_2)^{tr}$.  So it remains to count the $(1,2)$-good arrows in extreme position.  To do this, begin with a source monomial $t_s$ $=$ $x_2^u$, and construct all the vectors in extreme position that have source $t_s$ (refer to Figure \ref{fig:2vExtVec}).  One sees that the possible targets of such an arrow are the target monomials on the horizontal steps of the diagram that lie above and to the right of $t_s$, and a moment's reflection shows that the number of these is equal to the number of monomials in $\script{O}$ that lie on the same row as $t_s$.  It immediately follows from this that the number of $(1,2)$-good arrows in extreme position, and hence the number of spinal $(1,2)$-good tuples $(d_1, d_2)^{tr}$, is equal to $\mu$, the number of monomials in $\script{O}$, as the lemma asserts.   
\qedsymbol
\medskip

We now define a linear ordering $\prec$ on the arrows in extreme position.  First we list the monomials along the left edge of $\script{O}$ by ascending $x_2$-degree, as follows:
\[
  t_{p_0} = 1,\ t_{p_1} = x_2,\ t_{p_2} = x_2^2,\  \dots\ t_{p_h} = x_2^h.
\]
A $(1,2)$-good arrow in extreme position has the form $t_{p_s}$ $\rightarrow$ $m$, for some $1$ $\leq$ $s$ $\leq$ $h$ and $(1,2)$-target monomial $m$ as previously described.  The ordering is then defined as follows (see Figure \ref{fig:2vExtVec} for examples):
\[
   \left( t_{p_s} \rightarrow m \right)  \prec \left(t_{p_{s'}} \rightarrow m' \right) \Longleftrightarrow (s > s') \text{ or } ((s = s') \text{ and } (x_1\text{-deg}(m) < x_1\text{-deg}(m'))).
\]

\begin{lem} \label{lem:heart}
     Let $\rho^{1,2}_{p_s,q}$ be extreme, $\arr(\rho^{1,2}_{p_s,q})$ $=$ $\mathbf{a}$ $=$ $(t_{p_s}$ $\rightarrow$ $x_1 x_2 t_q)$ the associated extreme arrow with $(1,2)$-good direction vector $(d_1,d_2)^{tr}$, and $T$ $=$ $T_{\Pi,1}$ a trace syzygy with $\mdeg(\Pi)$ $=$ $(d_1,d_2)^{tr}$, so that, in particular, $\rho^{1,2}_{p_s,q}$ is in the spine of $T$.  Then for any other $(1,2)$-good vector $\mathbf{a'}$ $=$ $(t_{p_{s'}}$ $\rightarrow$ $x_1 x_2 t_{q'})$ in extreme position with $\mathbf{a}$ $\prec$ $\mathbf{a'}$, we have that the associated extreme $\rho^{1,2}_{p_{s'},q'}$ $=$ $\rho_{\iota'}$ will have component $\kappa_{\iota'}$ $=$ $0$ in $T$.  
\end{lem}

\emph{Proof:}
    If $\rho^{1,2}_{p_{s'},q'}$ is to appear nontrivially in $T$, then it must appear with a non-zero coefficient in the trace of one of the summands
\begin{equation} \label{eqn:summand}
    \gfmm_{l_1}\, \gfmm_{l_2}\, \dots\, \gfmm_{l_{v-1}}\, [\gfmm_1, \gfmm_{l_v}]\, \gfmm_{l_{v+1}}\, \dots\, \gfmm_{l_{s-1}}
\end{equation}
of (\ref{eqn:TSyzSum}), where each $l_i$ is either $1$ or $2$, and $l_v$ $=$ $2$.  However, by the proof of Proposition \ref{prop:Tspines}, we know that $\rho^{1,2}_{p_s,q}$ has coefficient $1$ in the trace of the matrix (\ref{eqn:summand}), and that this term comes from the $(r,r)$-th component of the product, where $r$ is the index of the monomial
\[
    t_r\ =\ (x_{l_1}\, x_{l_2}\, \dots\, x_{l_{v-1}}) t_{p_s}\ \Rightarrow\ (x_{l_{v+1}}\, x_{l_{v+2}}\, \dots\, x_{l_{s-1}}) t_r\ =\ t_q.
\]
We now consider what is required for $\rho^{1,2}_{p_{s'},q'}$ to appear with a non-zero coefficient in the trace of (\ref{eqn:summand}).  It would of course have to appear in the $(u,u)$-th component, which is given by 
\begin{equation}\label{eqn:uthComp}
  e_u \cdot  \left( \gfmm_{l_1}\, \gfmm_{l_2}\, \dots\, \gfmm_{l_{v-1}}\, [\gfmm_1, \gfmm_{l_v}]\, \gfmm_{l_{v+1}}\, \dots\, \gfmm_{l_{s-1}} \right) \cdot e_u^{tr},
\end{equation}
where $e_u$ denotes the $u$-th standard unit vector of length $\mu$.  
Proceeding as in the proof of Proposition \ref{prop:Tspines}, we suppose that
\[
    \left( \gfmm_{l_{v+1}}\dots \, \gfmm_{l_{s-1}} \right) \cdot e_u^{tr}\ =\ \mathbf{\gamma}^{tr} = (\gamma_1, \gamma_2, \dots, \gamma_{\mu})^{tr},\ \ \gamma_i \in \Zc. 
\]
Then $\rho^{1,2}_{p_{s'},q'}$ lies in the $p_{s'}$-th component (only) of
\[
    [\gfmm_1, \gfmm_2)] \cdot \mathbf{\gamma}^{tr}\ =\ (\rho^{1,2}_{ij})\cdot \mathbf{\gamma}^{tr},
\]
which is equal to $\sum_{j=1}^{\mu}\rho^{1,2}_{p_{s'},j}\cdot \gamma_j$.  To determine the coefficient of $\rho^{1,2}_{p_{s'},q'}$ in the trace, it remains to compute 
\[
  \begin{array}{cl}
   {} &  e_u \cdot \left( \gfmm_{l_1}\, \gfmm_{l_2}\dots \gfmm_{l_{v-1}} \right) \cdot (\sum_{j=1}^{\mu}\rho^{1,2}_{p_{s'},j}\cdot \gamma_j)e_{p_{s'}}^{tr}\vspace{.05in}\\
 = & (\sum_{j=1}^{\mu}\rho^{1,2}_{p_{s'},j}\cdot \gamma_j) \cdot e_u \cdot \left( \gfmm_{l_1}\, \gfmm_{l_2}\dots \gfmm_{l_{v-1}} \right) \cdot e_{p_{s'}}^{tr}.
  \end{array}
\]
But we are assuming that $\mathbf{a}$ $\prec$ $\mathbf{a'}$, which implies in particular that $s'$ $\leq$ $s$; that is, $t_{p_{s'}}$ either equals or lies directly below $t_{p_{s}}$ in the lattice of monomials.  It is then clear that $(x_{l_1} x_{l_2} \cdots x_{l_{v-1}})\cdot t_{p_{s'}}$ $=$ $t_{r'}$ $\in$ $\script{O}$, with $t_{r'}$ either equal to or lying directly below $t_r$ in the lattice.  The above product then reduces to 
\[
    (\sum_{j=1}^{\mu}\rho^{1,2}_{p_{s'},j}\cdot \gamma_j) \cdot e_u  \cdot e_{r'}^{tr},
\]
which equals $0$ unless $u$ $=$ $r'$.  In the case where $s'$ $=$ $s$ $\Rightarrow$ $r'$ $=$ $r$, we already know that  
\[
    \left( x_{l{v+1}}\, x_{l_{v+2}}\, \dots\, x_{l_{s-1}}\right) t_{r} \ =\ t_q. 
\]
It follows that only $\rho$'s of the form $\rho^{1,2}_{p_1,q}$, $1$ $\leq$ $p_1$ $\leq$ $\mu$, can appear in the $(u = r)$-th summand (\ref{eqn:uthComp}) of the trace; hence, $\rho^{1,2}_{p_{s'},q'}$ cannot appear, since $\mathbf{a}$ $\prec$ $\mathbf{a'}$ implies $q'$ $\neq$ $q$ in this case.  In the case where $s'$ $<$ $s$ and $t_{r'}$ lies directly below $t_r$ in the lattice of monomials, we must have that
\[
    \left(x_{l{v+1}}\, x_{l_{v+2}}\, \dots\, x_{l_{s-1}}\right) t_{r'}\ =\ t_{q''} \in \script{O}
\]
since the product is a monomial that lies directly below $t_q$ in the lattice.  One now sees that only $\rho$'s of the form $\rho^{1,2}_{p_1,q''}$ can appear in the summand (\ref{eqn:uthComp}).  If $\rho^{1,2}_{p_{s'},q'}$ is to appear, we must have $t_{q''}$ $=$ $t_{q'}$; however, it is then clear that $x_2 t_{q'}$ $\in$ $\script{O}$, which implies that $\mathbf{a'}$ $=$ $(t_{p_{s'}} \rightarrow x_1 x_2 t_{q'})$ is not extreme, contrary to hypothesis, since the target monomial $x_1 x_2 t_{q'}$ is not on a horizontal step of the diagram.  This completes the proof of the lemma. 
\qedsymbol
\medskip

Before harvesting our desired (and final) result, we recall that we are assuming that $\operatorname{char}(\gf)$ $=$ $0$ in this section, a hypothesis we have not yet used.

\begin{thm} \label{prop:compInt}
    Let $\script{O}$ $\subseteq$ $\mathbb{T}^2$ be an order ideal.  Then the ideal $I(\mathbb{B}_{\script{O}})$ of the border basis scheme 
\[
      \mathbb{B}_{\script{O}} = \Spec(\gf[\mathbf{c}]/I(\mathbb{B}_{\script{O}}))\ \subseteq\ \mathbb{A}^{\mu \nu} ,\ \ I(\mathbb{B}_{\script{O}}) = (\rho^{1,2}_{pq}),\ 1 \leq p, q \leq \mu,
\]
is generated by the $\rho^{1,2}_{pq}$ that are neither trivially $0$ nor extreme, and the number of these is $(\nu-2)\mu$.  Consequently, $\mathbb{B}_{\script{O}}$ is a complete intersection in $\Spec(\gf[\mathbf{c}])$.
\end{thm}

\emph{Proof:}
    By Lemma \ref{lem:numRhos}, $(\nu - 1)\mu$ of the $\rho^{1,2}_{pq}$ are not trivially $0$, and we know that these generate $I(\mathbb{B}_{\script{O}})$, by definition.  Furthermore, by Lemma \ref{lem:NumTSyz}, there are $\mu$ spinal multi-degrees $\mathbf{d}_{\ell}$ $=$ $(d_{\ell,1},d_{\ell,2})$, each of which corresponds to an arrow $\mathbf{a_{\ell}}$ $=$ $(t_{p_{s(\ell)}}$ $\rightarrow$ $m_{\ell})$ in extreme position, and an extreme element $\rho^{1,2}_{p_{s(\ell)},q_{\ell}}$ in the spine of $T_{\ell}$ $=$ $T_{\Pi_{\ell},1}$, where $\mdeg(\Pi_{\ell})$ $=$ $\mathbf{d}_{\ell}$, $1$ $\leq$ $\ell$ $\leq$ $\mu$; we assign the indices so that
\[
    \mathbf{a_{\mu}} \prec \mathbf{a_{\mu-1}} \prec \dots \prec \mathbf{a_{1}}.    
\]
Now Lemma \ref{lem:heart} tells us that $T_{\Pi_{\ell},1}$, which has the extreme $\rho^{1,2}_{p_{s(\ell)},q_{\ell}}$ in its spine, can involve no other extreme $\rho^{1,2}_{p_{s(\ell')},q_{\ell'}}$ with $\ell'$ $<$ $\ell$.  Moreover, since the coefficient of $\rho^{1,2}_{p_{s(\ell)},q_{\ell}}$ $=$ $\rho_{\iota_{\ell}}$ in $T_{\Pi_{\ell},1}$ is $\kappa_{\iota_{\ell}}$ $=$ $d_{\ell,2}$ $=$ (number of $2$'s in $\Pi_{\ell}$), and $\operatorname{char}(\gf)$ $=$ $0$, we can express $\rho^{1,2}_{p_{s(\ell)},q_{\ell}}$ as a $\gf[\mathbf{c}]$-linear combination of non-extreme $\rho$'s and extreme $\rho^{1,2}_{p_{s(\ell')},q_{\ell'}}$ with $\ell'$ $>$ $\ell$.  Therefore, starting with $\ell$ $=$ $\mu$ and working by descending induction, we obtain that every extreme $\rho^{1,2}_{p_{s(\ell)},q_{\ell}}$ lies in the ideal generated by the $\rho$'s that are neither trivially $0$ nor extreme, of which there are $(\nu-1)\mu$ $-$ $\mu$ $=$ $(\nu-2)\mu$; this completes the proof.
\qedsymbol

\begin{rem} \label{rem:sec9Eg}
    We have already seen two examples of Proposition \ref{prop:compInt}.  The first is a special case of Example \ref{exmp:zeroRhos}: 
\[
    \script{O}\ =\ \{ 1\} \subseteq \mathbb{T}^2,\ \ \partial \script{O} \ =\ \{x_1, x_2 \}.
\]
In this case, there is one $(1,2)$-target monomial $m$ $=$ $x_1 x_2$, one spinal multi-degree $(1,1)$, and one extreme arrow $(1$ $\rightarrow$ $m)$.  The associated trace syzygy is
\[
    T_{\left< 1,2 \right>,1} = \Tr \left( [\gfmm_1, \gfmm_2]\right) = \Tr\left(\rho^{1,2}_{1,1} \right) = \rho^{1,2}_{1,1}.
\]
Thus, in this case, there is only one $\rho$ that is both extreme and not trivially $0$, and the trace syzygy tells us that $\rho^{1,2}_{1,1}$ $=$ $0$, a result previously obtained by direct computation.  Hence, the ideal $I(\mathbb{B}_{\script{O}})$ is indeed generated by $(\nu-2)\mu$ $=$ $0 \cdot 1$ $=$ $0$ polynomials. 

 The second example occurs in Section \ref{sec:TSyzEx1}, where we considered the order ideal 
\[
    \script{O}\ =\  \{ 1,\, x_1,\, x_2 \} \subseteq \mathbb{T}^2.
\]
Indeed, (\ref{eqn:threeRhos}) states that the ideal $I(\mathbb{B}_{\script{O}})$ is generated by the three non-extreme $\rho$'s, and that the three extreme $\rho$'s can be omitted from the list of generators. Figure \ref{fig:1x1x2} shows the three extreme arrows (and one non-extreme arrow) in this case. 
 \end{rem}

\bibliographystyle{amsplain}
\bibliography{refs}

\end{document}